\documentstyle[11pt,amssymb]{article}

\begin{document}

\noindent
{\bf REPRESENTATIONS OF THE $q$-DEFORMED ALGEBRA
$U'_q({\rm so}_4)$}

\vskip 15 pt

{\sl M. Havl\'\i\v cek}

{\it Department of Mathematics, FNSPE, Czech Technical University}

{\it CZ-120 00, Prague 2, Czech Republic}
\medskip

{\sl A. U. Klimyk}

{\it Institute for Theoretical Physics, Kiev 252143, Ukraine}
\medskip

{\sl S. Po\v sta}

{\it Department of Mathematics, FNSPE, Czech Technical University}

{\it CZ-120 00, Prague 2, Czech Republic}

\vskip 30 pt

\begin{abstract}
We study the nonstandard $q$-deformation $U'_q({\rm so}_4)$ of the
universal enveloping algebra $U({\rm so}_4)$ obtained by deforming the
defining relations for skew-symmetric generators of $U({\rm so}_4)$. This
algebra is used in quantum gravity and algebraic topology. We construct
a homomorphism $\phi$ of $U'_q({\rm so}_4)$ to the certain nontrivial
extension of the Drinfeld--Jimbo quantum algebra
$U_q({\rm sl}_2)^{\otimes 2}$ and show that this homomorphism is an
isomorphism. By using this homomorphism we construct irreducible
finite dimensional representations of the classical type and of the
nonclassical type for the algebra $U'_q({\rm so}_4)$. It is proved that
for $q$ not a root of unity each irreducible finite dimensional
representation
of $U'_q({\rm so}_4)$ is equivalent to one of these representations. We
prove
that every finite dimensional representation of $U'_q({\rm so}_4)$ for
$q$ not a root of unity is completely reducible. It is shown how to
construct (by using the homomorphism $\phi$)
tensor products of irreducible representations of $U'_q({\rm so}_4)$.
(Note that no Hopf algebra structure is known for $U'_q({\rm so}_4)$.)
These tensor products are decomposed into irreducible constituents.
\end{abstract}

\noindent
{\sf 1. INTRODUCTION}
\medskip

In [1] it was introduced a $q$-deformation of the universal enveloping
algebra $U({\rm so}_n)$ which uses a realization of the Lie algebra
${\rm so}_n$ by skew-symmetric matrices (rather than by root elements
and the diagonal Cartan subalgebra). If $I_{ij}=E_{ij}-E_{ji}$ are the
skew-symmetric matrices of ${\rm so}_n$, where $(E_{ij})_{sr}:=
\delta _{is}\delta _{jr}$, then ${\rm so}_n$ is generated by the matrices
$I_{21}, I_{32},\cdots ,I_{n,n-1}$ (other basis matrices are obtained by
taking commutators of these matrices). The following Serre type theorem
is true (see [2]): {\it The universal enveloping algebra
$U({\rm so}_n)$ is isomorphic to the associative algebra generated by
(abstract) elements $I_{21}, I_{32},\cdots ,I_{n,n-1}$ satisfying the
relations}
$$
I_{i,i-1}^2I_{i+1,i}-2 I_{i,i-1}I_{i+1,i}I_{i,i-1} +I_{i+1,i}
I^2_{i,i-1} =-I_{i+1,i} ,
$$        $$
I_{i,i-1}I^2_{i+1,i}-2I_{i+1,i}I_{i,i-1}I_{i+1,i} +I^2_{i+1,i}
I_{i,i-1} =-I_{i,i-1} ,
$$       $$
I_{i,i-1}I_{j,j-1}= I_{j,j-1}I_{i,i-1}\ \ \ \ {\rm for}\ \ \ \ |i-j|>1.
$$
Now we $q$-deform these relations by $2\to [2]_q\equiv q+q^{-1}$.
As a result, we obtain the associative algebra generated by elements
$I_{21}, I_{32},\cdots ,I_{n,n-1}$ satisfying the relations
$$
I_{i,i-1}^2I_{i+1,i}-(q+q^{-1}) I_{i,i-1}I_{i+1,i}I_{i,i-1}+I_{i+1,i}
I^2_{i,i-1} =-I_{i+1,i} ,
$$        $$
I_{i,i-1}I^2_{i+1,i}-(q+q^{-1})I_{i+1,i}I_{i,i-1}I_{i+1,i}+I^2_{i+1,i}
I_{i,i-1} =-I_{i,i-1} ,
$$       $$
I_{i,i-1}I_{j,j-1}= I_{j,j-1}I_{i,i-1}\ \ \ \ {\rm for} \ \ \ \
|i-j|>1.
$$
We denote this algebra by $U'_q({\rm so}_n)$. Associative algebras
isomorphic to $U'_q({\rm so}_n)$ appear in quantum gravity [3-5],
in discrete Schr\"odinger equation [6], in algebraic topology [7-8],
in the theory of $q$-orthogonal polynomials [9] and in the theory of
$q$-Laplace operators and $q$-harmonic polynomials [10]. For this reason,
studying the algebra $U'_q({\rm so}_n)$ (especially for small numbers
$n$) is of great importance. There are several problems which have to
be solved. The most important are the following ones:
\medskip

(a) relation of $U'_q({\rm so}_n)$ to Drinfeld--Jimbo quantum algebras;

(b) structure of the algebra $U'_q({\rm so}_n)$ (including
explicit form of the center, Casimir elements,  automorphism group, etc);

(c) construction and classification of irreducible finite dimensional
representations.
\medskip

\noindent
Part of the main problems are solved. For example, it is shown (see
[11] and [12]) that $U'_q({\rm so}_n)$ can be embedded as a subalgebra
into the Drinfeld--Jimbo quantum algebra $U_q({\rm sl}_n)$. (Remind that the
Drinfeld--Jimbo quantum algebra $U_q({\rm so}_n)$ is not contained in
$U_q({\rm sl}_n)$). The main classes of irreducible finite dimensional
representations of $U'_q({\rm so}_n)$ are constructed for $q$ not a root
of unity (see [13] and [14]) and for $q$ a root of unity (see
[12] and [15]). But
the corresponding classification theorem for the representation theory
of $U'_q({\rm so}_n)$ is not proved. Casimir elements were constructed
for $q$ not a root of unity (see [16] and [17]) and for $q$ a root of
unity (see [18]). But it is not known whether they generate the center.

Most problems are solved [19-22] for the algebra $U'_q({\rm so}_3)$.
In particular, it was shown that $U'_q({\rm so}_3)$ can be
embedded into a certain extension of the Drinfeld--Jimbo algebra
$U_q({\rm sl}_2)$ (but $U'_q({\rm so}_3)$ is not isomorphic to
$U_q({\rm sl}_2)$). The classification theorem for representation theory
of $U'_q({\rm so}_3)$ were proved [22]. It was shown that
the automorphism group of $U'_q({\rm so}_3)$ contains a group isomorphic
to the modular group $SL(2,{\Bbb Z})$.

The aim of this paper is to solve the main problems for the algebra
$U'_q({\rm so}_4)$ naturally appearing in algebraic topology [8].
(Note that unlike to the classical case, the algebra $U'_q({\rm so}_4)$
is not isomorphic to the product of two copies of $U'_q({\rm so}_3)$).
The algebra $U'_q({\rm so}_4)$ and its irreducible finite
dimensional representations were studied in several papers [23-25].
Nevertheless, the main problems were not solved.
The main results of this paper are the following:
\medskip

(I) We construct a homomorphism from $U'_q({\rm so}_4)$ to some
extension of the Drinfeld--Jimbo quantum algebra
$U_q({\rm sl}_2)^{\otimes 2}\equiv U_q({\rm sl}_2)\otimes
U_q({\rm sl}_2)$. It is shown that this homomorphism is injective.
Thus, $U'_q({\rm so}_4)$ is embedded into this extension of
$U_q({\rm sl}_2)^{\otimes 2}$. This solves for $U'_q({\rm so}_4)$
the problem (a).

(II) It is proved the theorem classifying irreducible finite dimensional
representations of $U'_q({\rm so}_4)$, when $q$ is not a root of unity.
According to this theorem, irreducible representations of the classical
type ($q$-analogue of the irreducible representations of the Lie algebra
${\rm so}_4$) and irreducible representations of the nonclassical type
(these representations do not have any classical analogue) exhaust all
irreducible finite dimensional representations of $U'_q({\rm so}_4)$.
This solves for $U'_q({\rm so}_4)$, when $q$ is not a root of unity,
the problem (c).

(III) It is proved that if $q$ is not a root of unity, then any finite
dimensional representation of $U'_q({\rm so}_4)$ is completely reducible.

(IV) It is shown how to construct tensor products of irreducible
finite dimensional representations of $U'_q({\rm so}_4)$.
(Note that no Hopf algebra structure is known on $U'_q({\rm so}_4)$.)
\medskip

In sections 2--4 $q$ is any complex number different from $\pm 1$.
In other sections it is assumed that $q$ is not a root of unity.
\medskip

\noindent
{\sf 2. THE ALGEBRA $U'_q({\rm so}_4)$}
\medskip

We first define the $q$-deformed algebra $U'_q({\rm so}_3)$
which is a subalgebra of $U'_q({\rm so}_4)$. The algebra
$U'_q({\rm so}_3)$ is obtained [26] by a $q$-deformation of the standard
commutation relations
$$
[I_{21},I_{32}]=I_{31},\ \ \ \
[I_{32},I_{31}]=I_{21},\ \ \ \
[I_{31},I_{21}]=I_{32}
$$
of the Lie algebra ${\rm so}_3$. So, the algebra $U'_q({\rm so}_3)$
is the complex associative algebra (with unit element) generated by
the elements $I_{21}$, $I_{32}$, $I_{31}$ satisfying the defining
relations
$$
[I_{21},I_{32}]_q: = q^{1/2}I_{21}I_{32}-q^{-1/2}I_{32}I_{21}
=I_{31}, \eqno (1)
$$    $$
[I_{32},I_{31}]_q: = q^{1/2}I_{32}I_{31}-q^{-1/2}I_{31}I_{32}
=I_{21}, \eqno (2)
$$      $$
[I_{31},I_{21}]_q: = q^{1/2}I_{31}I_{21}-q^{-1/2}I_{21}I_{31}
=I_{32}. \eqno (3)
$$
Note that by (1) the element $I_{31}$ is not independent: it depends
on the elements $I_{21}$ and $I_{32}$. Substituting the expression (1)
for $I_{31}$ into (2) and (3) we obtain the relations
$$
I_{21}I^2_{32} -(q+q^{-1})I_{32}I_{21}I_{32}+I^2_{32}I_{21}=-I_{21},
 \eqno (4)
$$       $$
I_{32}I^2_{21} -(q+q^{-1})I_{21}I_{32}I_{21}+I^2_{21}I_{32}=-I_{32} .
 \eqno (5)
$$
The relations (4) and (5) restore the relations (2) and (3) if to
introduce the element $I_{31}$ defined by (1).
The algebra $U'_q({\rm so}_3)$ can be defined as the associative
algebra generated by the elements $I_{21}$ and $I_{32}$ satisfying
the defining relations (4) and (5).

Starting from the definition of $U'_q({\rm so}_3)$ by relations (4)
and (5), we give the following definition of the $q$-deformed
algebra $U'_q({\rm so}_4)$. It is an associative algebra (with unit
element) generated by the elements $I_{21}$, $I_{32}$, $I_{43}$
satisfying the defining relations
$$
I_{21}I^2_{32} -(q+q^{-1})I_{32}I_{21}I_{32}+I^2_{32}I_{21}=-I_{21},
 \eqno (6)
$$       $$
I_{32}I^2_{21} -(q+q^{-1})I_{21}I_{32}I_{21}+I^2_{21}I_{32}=-I_{32} ,
 \eqno (7)
$$    $$
I_{32}I^2_{43} -(q+q^{-1})I_{43}I_{32}I_{43}+I^2_{43}I_{32}=-I_{32},
 \eqno (8)
$$       $$
I_{43}I^2_{32} -(q+q^{-1})I_{32}I_{43}I_{32}+I^2_{32}I_{43}=-I_{43} ,
 \eqno (9)
$$     $$
I_{21}I_{43}-I_{43}I_{21}=0. \eqno (10)
$$
It is clear that $U'_q({\rm so}_4)$ contains at least
two subalgebras isomorphic to $U'_q({\rm so}_3)$. The first one is
generated by $I_{21}$ and $I_{32}$, and the second one by $I_{32}$ and
$I_{43}$.

We can introduce in $U'_q({\rm so}_4)$ the elements $I_{31}$, $I_{42}$ and
$I_{41}$. They are defined as in (1):
$$
I_{31}:=[I_{21},I_{32}]_q,\ \ \ \ \ I_{42}:= [I_{32},I_{43}]_q,
\ \ \ \ \  I_{41}:=[I_{21},I_{42}]_q=[I_{31},I_{43}]_q,
$$
where $[A,B]_q:=q^{1/2}AB-q^{-1/2}BA$ is the $q$-commutator of $A$ and $B$.
Then the elements $I_{ij}$, $4\ge i>j\ge 1$, satisfy the relations [13]
$$
[I_{21},I_{32}]_q=I_{31},\ \ \  [I_{32},I_{31}]_q=I_{21},\ \ \
[I_{31},I_{21}]_q=I_{32},   \eqno (11)
$$   $$
[I_{32},I_{43}]_q=I_{42},\ \ \  [I_{43},I_{42}]_q=I_{32},\ \ \
[I_{42},I_{32}]_q=I_{43},    \eqno (12)
$$   $$
[I_{31},I_{43}]_q=I_{41},\ \ \  [I_{43},I_{41}]_q=I_{31},\ \ \
[I_{41},I_{31}]_q=I_{43},     \eqno (13)
$$     $$
[I_{21},I_{42}]_q=I_{41},\ \ \  [I_{42},I_{41}]_q=I_{21},\ \ \
[I_{41},I_{21}]_q=I_{42},      \eqno (14)
$$     $$
[I_{21},I_{43}]=0,\ \ \  [I_{32},I_{41}]=0,\ \ \
[I_{42},I_{31}]=(q-q^{-1})(I_{21}I_{43}-I_{32}I_{41}).  \eqno (15)
$$
As in the case of the algebra $U'_q({\rm so}_3)$, the relations (11)--(15)
are equivalent to the relations (6)--(10). Note that the relations
(11)--(15) define the algebra appearing in algebraic topology [8].

Four sets of relations (11)--(14) gives four subalgebras
of $U'_q({\rm so}_4)$ isomorphic to $U'_q({\rm so}_3)$. They are
generated by triples
$$
(I_{21},I_{32},I_{31}),\ \ \ \
(I_{32},I_{43},I_{42}),\ \ \ \
(I_{31},I_{43},I_{41}),\ \ \ \
(I_{21},I_{42},I_{41}),
$$
respectively.

The Poincar{\'e}--Birkhoff--Witt theorem is true for $U'_q({\rm so}_4)$.
It can be formulated as: {\it The elements
$$
I^{m_{31}}_{31}I^{m_{32}}_{32}I^{m_{41}}_{41}
I^{m_{42}}_{42}I^{m_{43}}_{43}I^{m_{21}}_{21},\ \ \ \
m_{ij}=0,1,2,\cdots ,
$$
form a basis of $U'_q({\rm so}_4)$.} This theorem is proved by means
of the diamond lemma (see [27], subsection 4.1.5). As in the case of
ordinary simple Lie algebra, the same theorem holds for any other
ordering of the six generators.

We shall need Casimir elements of $U'_q({\rm so}_4)$. In order to give
them we introduce also the elements
$$
I^-_{31}:=[I_{21},I_{32}]_{q^{-1}}, \ \ \ \
I^-_{42}:= [I_{32},I_{43}]_{q^{-1}},\ \ \ \
I^-_{41}:=[I_{31},I_{43}]_{q^{-1}},
$$
where $[A,B]_{q^{-1}}:=q^{-1/2}AB-q^{1/2}BA$. Then ( see [16])
$$
C_4 =q^{-1}I_{21}I_{43}- I_{31} I_{42} + q I_{32}I_{41}
$$       $$
C'_4 =q^{-2}I_{21}^2+I^2_{32} +q^2 I_{43}^2 +
q^{-1}I_{31}I^-_{31} + q I_{42}I^-_{42} +I_{41}I^-_{41} ,
$$
are two independent elements of the center of the algebra $U'_q({\rm
so}_4)$.
Using the Poincar{\'e}--Birkhoff--Witt theorem the element $C'_4$ can
be represented in the form
$$
C'_4=q^2 I_{21}^2+I^2_{41}+I^2_{32}+q^{-2}(I^2_{43}+I^2_{21}+I^2_{31})
-(q-q^{-1})q^{-3/2}(I_{31}I_{32}I_{21}+I_{31}I_{41}I_{43})-
$$  $$
-(q-q^{-1})q^{1/2} (I_{32}I_{42}I_{43}+I_{41}I_{42}I_{21})+(q-q^{-1})^2
I_{32}I_{41}I_{43}I_{21} .
$$

\noindent
{\sf 3. THE ALGEBRA $U_q({\rm sl}_2)^{\otimes 2,{\rm ext}}$}
\medskip

Let $e_1,f_1,q^{H_1}$ and $e_2,f_2,q^{H_2}$ be generating elements of two
copies of the quantum algebra $U_q({\rm sl}_2)$ satisfying the relations
$$
q^{H_i}e_i=qe_iq^{H_i},\ \ \ \
q^{H_i}f_i=q^{-1}f_iq^{H_i},\ \ \ \
[e_i,f_i]=\frac{q^{2H_i}-q^{-2H_i}}{q-q^{-1}} .
$$
The expressions
$$
c_i =e_if_i + \frac{q^{2H_i-1}+q^{-2H_i+1}}{(q-q^{-1})^2},\ \ \ \ \
i=1,\, 2,
$$
give Casimir elements of these algebras $U_q({\rm sl}_2)$.
The comultiplication $\Delta$ is given in $U_q({\rm sl}_2)$ by the
formulas
$$
\Delta (q^{\pm H_i})=q^{\pm H_i}\otimes q^{\pm H_i},\ \ \ \
\Delta (e_i)=e_i\otimes q^{ H_i}+
q^{- H_i}\otimes e_i,\ \ \ \
\Delta (f_i)=f_i\otimes q^{H_i}+
q^{- H_i}\otimes f_i .
$$

Let us consider the polynomials
$$
p_i(x_i)=q^{-1}x^4_i-c_i(q-q^{-1})^2x^2_i+q ,
$$
where $c_i$ are the Casimir elements.
They are irreducible in $U_q({\rm sl}_2)$, that is, there exists
no element $a\in U_q({\rm sl}_2)$ such that $p_i(a)=0$.
Therefore, we can define the quadric algebraic extension
${\hat U}_q({\rm sl}_2)$ of the algebra $U_q({\rm sl}_2)$ by means of
the element $x_i$ commuting with all elements of $U_q({\rm sl}_2)$:
$$
{\hat U}_q({\rm sl}_2)=\{ a_3x_i^3 +a_2 x_i^2 + a_1x_i+a_0\, |\,
a_j\in U_q({\rm sl}_2)\} ,
$$
assuming that $p_i(x_i)=0$, that is, $x_i^4=qc_i(q-q^{-1})^2x_i^2-q^2$.
This equation is equivalent to the following one
$$
c_i=\frac{x_i^2q^{-1}+x_i^{-2}q}{(q-q^{-1})^2} . \eqno (16)
$$
Note that the element $x_i$ has an inverse in ${\hat U}_q({\rm sl}_2)$
since
$$
x_i(-x_i^3q^{-1}+c_i(q-q^{-1})^2x_i)q^{-1}=1,
$$
that is,
$$
x_i^{-1}=(-x_i^3q^{-1}+c_i(q-q^{-1})^2x_i)q^{-1}.    \eqno (17)
$$

We consider two algebras ${\hat U}_q({\rm sl}_2)$ generated by the
elements $e_1,f_1,q^{\pm H_1}, x_1$ and
by the elements $e_2,f_2,q^{\pm H_2}, x_2$, respectively.
Let ${\hat U}_q({\rm sl}_2)^{\otimes 2}$ be the tensor product of
these algebras. Then we extend (in the sense of [28]) this algebra
${\hat U}_q({\rm sl}_2)^{\otimes 2}$ by adding to it the
commuting elements
$$
(q^{H_1}q^{H_2}q^j+q^{-H_1}q^{-H_2}q^{-j})^{-1},\ \ \ \
(q^{H_1}q^{-H_2}q^j+q^{-H_1}q^{H_2}q^{-j})^{-1},\ \ \ \
j=0,\pm 1,\pm 2,\cdots .
\eqno (18)
$$
This extended algebra will be denoted by
${\hat U}_q({\rm sl}_2)^{\otimes 2, {\rm ext}}$. It is the
associative algebra generated by the elements $e_i,f_i,q^{\pm H_i}, x_i$,
$i=1,\, 2$, and by elements (18) such that
$e_1,f_1,q^{\pm H_1}, x_1$ and
$e_2,f_2,q^{\pm H_2}, x_2$ satisfy the relations determined in the
algebra ${\hat U}_q({\rm sl}_2)$, each of the elements
$e_1,f_1,q^{\pm H_1}, x_1$ commute with each of the elements
$e_2,f_2,q^{\pm H_2}, x_2$, each of the elements $q^{\pm H_1}$ and
$q^{\pm H_2}$ commutes with each of elements (18), and
$$
(q^{H_1+H_2+j}+q^{-H_1-H_2-j})^{-1}e_i=
e_i(q^{H_1+H_2+j+1}+q^{-H_1-H_2-j-1})^{-1},
$$       $$
(q^{H_1+H_2+j}+q^{-H_1-H_2-j})^{-1}f_i=
f_i(q^{H_1+H_2+j-1}+q^{-H_1-H_2-j+1})^{-1},
$$       $$
(q^{H_1-H_2+j}+q^{-H_1+H_2-j})^{-1}e_i=
e_i(q^{H_1-H_2+j+\varepsilon }+q^{-H_1+H_2-j-\varepsilon })^{-1},
$$       $$
(q^{H_1-H_2+j}+q^{-H_1+H_2-j})^{-1}f_i=
f_i(q^{H_1-H_2+j-\varepsilon }+q^{-H_1+H_2-j+\varepsilon })^{-1},
$$       $$
(q^{H_1+H_2+j}+q^{-H_1-H_2-j})^{-1}
(q^{H_1+H_2+j}+q^{-H_1-H_2-j})=1,
$$       $$
(q^{H_1-H_2+j}+q^{-H_1+H_2-j})^{-1}
(q^{H_1-H_2+j}+q^{-H_1+H_2-j})=1,
$$
where $\varepsilon =1$ if $i=1$ and $\varepsilon =-1$ if $i=2$.

Let us find irreducible finite dimensional representations of the
algebra ${\hat U}_q({\rm sl}_2)^{\otimes 2,{\rm ext}}$ for $q$ not
a root of unity. For these values of $q$
the algebra $U_q({\rm sl}_2)$ has finite dimensional
irreducible representations $T_l\equiv T^{(1)}_l$, $T^{(-1)}_l$,
$T^{({\rm i})}_l$, $T^{(-{\rm i})}_l$, $l=0,\frac 12 ,1,\frac 32 ,
\cdots$, acting on the vector space ${\cal H}_l$ with basis $|l,m
\rangle$, $m=-l,-l+1,\cdots ,l$. These representations are given by
the formulas
$$
T^{(1)}_l(q^H) |l,m\rangle = q^m |l,m\rangle ,\ \ \ \
T^{(1)}_l(e) |l,m\rangle = [l-m] |l,m+1\rangle ,  \eqno (19)
$$      $$
T^{(1)}_l(f) |l,m\rangle = [l+m] |l,m-1\rangle ,   \eqno (20)
$$
where numbers in square brackets mean $q$-numbers determined by
$$
[a]=\frac{q^q-q^{-a}}{q-q^{-1}} ,
$$
and by the formulas
$$
T^{(-1)}_l(q^H) |l,m\rangle =-q^H |l,m\rangle ,\ \ \ \
T^{(-1)}_l(X)=  T^{(1)}_l(X) ,\ \ \ X=e,f,  \eqno (21)
$$      $$
T^{(\rm i)}_l(q^H) |l,m\rangle ={\rm i}q^H |l,m\rangle ,\ \ \ \
T^{(\rm i)}_l(e)=  T^{(1)}_l(e) ,\ \ \ \
T^{(\rm i)}_l(f)=  -T^{(1)}_l(f) ,    \eqno (22)
$$      $$
T^{(-{\rm i})}_l(q^H) |l,m\rangle =-{\rm i}q^H |l,m\rangle ,\ \ \ \
T^{(-{\rm i})}_l(e)=  T^{(1)}_l(e) , \ \ \ \
T^{(-{\rm i})}_f(f)=  -T^{(1)}_l(f)    \eqno (23)
$$
(see, for example, [27], chapter 3). The representations
$T^{(1)}_l$, $T^{(-1)}_l$, $T^{({\rm i})}_l$, $T^{(-{\rm i})}_l$,
$l=0,\frac 12 ,1,\frac 32 ,\cdots$, are pairwise nonequivalent and any
irreducible finite dimensional representation of $U_q({\rm sl}_2)$
is equivalent to one of these representations.
Values of the Casimir element $c$ on these representations are given by
$$
T^{(1)}_l(c)=T^{(-1)}_l(c)=\frac{q^{2l+1}+q^{-2l-1}}{(q-q^{-1})^2},
\ \ \ \ \ \ \
T^{({\rm i})}_l(c)=T^{(-{\rm i})}_l(c)=
-\frac{q^{2l+1}+q^{-2l-1}}{(q-q^{-1})^2}.
$$

Since the Casimir element of $U_q({\rm sl}_2)$ is multiple to the
unit operator on the space ${\cal H}_l$, then each of the
representations
$T^{(1)}_l$, $T^{(-1)}_l$, $T^{({\rm i})}_l$, $T^{(-{\rm i})}_l$,
$l=0,\frac 12 ,1,\frac 32 ,\cdots$, can be extended to a representation
of ${\hat U}_q({\rm sl}_2)$. In order to determine these extensions we have
to determine the operators $T^{(\varepsilon )}_l (x)$, $\varepsilon =\pm 1,
\pm {\rm i}$, corresponding to the element $x$ from (16).
It follows from (16) that
$$
T^{(\varepsilon )}_l (c) =\frac{ T^{(\varepsilon )}_l (x)^2q^{-1} +
T^{(\varepsilon )}_l (x)^{-2}q}{(q-q^{-1})^2} .
$$
If some operator $T^{(\varepsilon )}_l (x)$ is a solution of this equation,
then the operators
$$
{\tilde T}^{(\varepsilon )}_l (x)=-  T^{(\varepsilon )}_l (x),\ \
T^{(\varepsilon )}_l (x)^{-1}q,\ \  -T^{(\varepsilon )}_l (x)^{-1}q
$$
are also its solutions. Each of them can be taken for extension of the
representation $T^{(\varepsilon )}_l$ of $U_q({\rm sl}_2)$. Since the
element $x$ commute with all elements of $U_q({\rm sl}_2)$, then different
extensions of $T^{(\varepsilon )}_l$ (obtained by using different
solutions of the above equation) do not essentially differ from each other.
For this reason, we shall use only the solution
$T^{(\varepsilon )}_l(x)=q^{-l}I$, where $I$ is the unit operator.
We denote the extended representations of
${\hat U}_q({\rm sl}_2)$, extended by using the solution
$T^{(\varepsilon )}_l(x)=q^{-l}I$, by the
same symbols $T^{(1)}_l$, $T^{(-1)}_l$, $T^{({\rm i})}_l$,
$T^{(-{\rm i})}_l$.

It is clear that irreducible finite dimensional
representations of ${\hat U}_q({\rm sl}_2)^{\otimes 2}$ are
equivalent to the following ones:
$$
T_l^{(\varepsilon )}\otimes T_{l'}^{(\varepsilon ')},\ \ \ \
l,l'=0,\frac 12 ,1,\frac 32 ,\cdots ,\ \ \ \
\varepsilon ,\varepsilon '=1,-1,{\rm i},-{\rm i} ,
$$
where $T_l^{(\varepsilon )}$ and $T_{l'}^{(\varepsilon ')}$ are irreducible
representations of two copies of the algebra ${\hat U}_q({\rm sl}_2)$,
respectively. Now we wish to extend these representations of
${\hat U}_q({\rm sl}_2)^{\otimes 2}$ to representations of the algebra
${\hat U}_q({\rm sl}_2)^{\otimes 2,{\rm ext}}$ by using the relation
$$
(T_l^{(\varepsilon )}\otimes T_{l'}^{(\varepsilon ')})\left( (q^{H_1+H_2+j}+
q^{-H_1-H_2-j})^{-1}\right) =
$$    $$
=\left( q^jT_l^{(\varepsilon )}(q^{H_1})\otimes
T_{l'}^{(\varepsilon ')}(q^{H_2})+ q^{-j}
T_l^{(\varepsilon )} (q^{-H_1}) \otimes
T_{l'}^{(\varepsilon ')}(q^{-H_2})\right) ^{-1} ,
$$      $$
(T_l^{(\varepsilon )}\otimes T_{l'}^{(\varepsilon ')})\left( (q^{H_1-H_2+j}+
q^{-H_1+H_2-j})^{-1}\right) =
$$          $$
=\left( q^jT_l^{(\varepsilon )}(q^{H_1})\otimes
T_{l'}^{(\varepsilon ')}(q^{-H_2})+ q^{-j}
T_l^{(\varepsilon )} (q^{-H_1}) \otimes
T_{l'}^{(\varepsilon ')}(q^{H_2})\right) ^{-1} .
$$
Clearly, only those irreducible representations
$T_l^{(\varepsilon )}\otimes T_{l'}^{(\varepsilon ')}$ of
${\hat U}_q({\rm sl}_2)^{\otimes 2}$ can be extended
to representations of
${\hat U}_q({\rm sl}_2)^{\otimes 2, {\rm ext}}$
for which all the operators
$$
q^jT_l^{(\varepsilon )}(q^{H_1})\otimes
T_{l'}^{(\varepsilon ')}(q^{H_2})+ q^{-j}
T_l^{(\varepsilon )} (q^{-H_1}) \otimes
T_{l'}^{(\varepsilon ')}(q^{-H_2}) , \ \ \ \ j=0,\pm 1,\cdots ,
$$    $$
 q^jT_l^{(\varepsilon )}(q^{H_1})\otimes
T_{l'}^{(\varepsilon ')}(q^{-H_2})+ q^{-j}
T_l^{(\varepsilon )} (q^{-H_1}) \otimes
T_{l'}^{(\varepsilon ')})(q^{H_2}) , \ \ \ \ j=0,\pm 1,\cdots ,
$$
are invertible. From formulas (19)--(23) it follows that these operators
are always invertible for the representations $T_l^{(\varepsilon )}\otimes
T_{l'}^{(\varepsilon ')}$ such that $\varepsilon ,\varepsilon '=1,-1$ or
$\varepsilon ,\varepsilon ' ={\rm i},-{\rm i}$, and also for all the
representations $T_l^{(\varepsilon )}\otimes T_{l'}^{(\varepsilon ')}$
such that $l+l'$ is half-integral (but not integral) number and
$\varepsilon =\pm 1$, $\varepsilon '=\pm {\rm i}$ or
$\varepsilon =\pm {\rm i}$, $\varepsilon '=\pm1$. For the representations
$T_l^{(\varepsilon )}\otimes T_{l'}^{(\varepsilon ')}$ with $\varepsilon
=\pm 1$,
$\varepsilon ' =\pm {\rm i}$ or $\varepsilon =\pm {\rm i}$, $\varepsilon
'=\pm 1$
and $l+l'\in {\Bbb Z}$ some of these operators are not invertible since
they have zero eigenvalue. Denoting the extended representations by the
same symbols, we can formulate the following assertion:
\medskip

{\bf Theorem 1.} {\it If $q$ is not a root of unity, then the algebra}
${\hat U}_q({\rm sl}_2)^{\otimes 2,
{\rm ext}}$ {\it has irreducible finite dimensional representations}
$$
T_l^{(\varepsilon )}\otimes T_{l'}^{(\varepsilon ')} ,\ \ \ \ l,l'=0,\frac
12
,1,\frac 32 ,\cdots ,\ \ \
\varepsilon ,\varepsilon '=\pm 1 \ \ {\rm or}\ \ \varepsilon ,\varepsilon
'=\pm {\rm i},
$$     $$
T_l^{(\varepsilon )}\otimes T_{l'}^{(\varepsilon ')}, \ \ \ \ l+l'\in {\frac
12}{\Bbb Z},\ \ \
l+l'\notin {\Bbb Z},\ \ \ \varepsilon =\pm 1,\ \varepsilon '=\pm {\rm i}\ \
\
{\rm or}\ \ \ \varepsilon  =\pm {\rm i},\ \varepsilon '=\pm 1
$$
{\it (all four combinations of signs are possible).
Up to values of the operators corresponding to the elements
$x_1^k\otimes x_2^s$, $k,s=1,2,3$,
any irreducible finite dimensional representation of ${\hat U}_q({\rm
sl}_2)^{\otimes 2, {\rm ext}}$
is equivalent to one of these representations.}
\bigskip

\noindent
{\sf 4. THE ALGEBRA HOMOMORPHISM $U'_q({\rm so}_4)\to
{\hat U}_q({\rm sl}_2)^{\otimes 2, {\rm ext}}$}
\medskip

The aim of this section is to give in an explicit form the
algebra homomorphism of
$U'_q({\rm so}_4)$ to ${\hat U}_q({\rm sl}_2)^{\otimes 2,
{\rm ext}}$. This homomorphism is given by the following theorem:
\medskip

{\bf Theorem 2.} {\it There exists a unique algebra homomorphism}
$\phi : U'_q({\rm so}_4)\to {\hat U}_q({\rm sl}_2)^{\otimes 2, {\rm ext}}$
{\it such that}
$$
\phi (I_{21}) = {\rm i}[H_1+H_2]_q,\ \ \ \
\phi (I_{43}) = {\rm i}[H_1-H_2]_q, \eqno (24)
$$      $$
\phi (I_{32}) =-\frac{x_1^{-1}q^{-H_2+1}+x_1q^{H_2-1}}
{(q^{H_1+H_2-1}+q^{-H_1-H_2+1})(q^{H_1-H_2+1}+q^{-H_1+H_2-1})} E_2 +
$$      $$
+\frac{x_1^{-1}q^{H_2+1}+x_1q^{-H_2-1}}
{(q^{H_1+H_2+1}+q^{-H_1-H_2-1})(q^{H_1-H_2-1}+q^{-H_1+H_2+1})} F_2 +
$$     $$
+\frac{x_2^{-1}q^{-H_1+1}+x_2q^{H_1-1}}
{(q^{H_1+H_2-1}+q^{-H_1-H_2+1})(q^{H_1-H_2-1}+q^{-H_1+H_2+1})} E_1 -
$$      $$
-\frac{x_2^{-1}q^{H_1+1}+x_2q^{-H_1-1}}
{(q^{H_1+H_2+1}+q^{-H_1-H_2-1})(q^{H_1-H_2+1}+q^{-H_1+H_2-1})} F_1 .
\eqno (25)
$$

{\sl Proof.} We have to show that three elements $\phi (I_{21})$,
$\phi (I_{32})$ and $\phi (I_{43})$ from (24)  and (25)
satisfy the defining relations (6)--(10). It is made by direct
verification. Namely, we substitute the expressions (24) and (25)
for $\phi (I_{21})$, $\phi (I_{43})$, $\phi (I_{32})$ into
(6)--(10) and then permute the generating elements $(q^{H_i})^{\pm 1}$,
$e_i$, $f_i$ in numerators (using the defining relations of the
algebra $U_q({\rm sl}_2)$) reducing them to the form
$$
(q^{H_1})^r (q^{H_2})^s e_1^{a_1} e_2^{b_1} f_1^{a_2} f_2^{b_2} ,\ \ \ \ \
r,s\in {\Bbb Z},\ \ \ \ a_1,b_1,a_2,b_2 \in {\Bbb Z}_+.
$$
Then it is directly seen that the relations (7), (8) and (10) are
fulfilled. So, we have to prove the relations (6) and (9).
We cancel in these relations separately terms
ending with $e_1^2$, $e_2^2$, $f_1^2$ and $f_2^2$.
Now in the relation (6) we cancel terms ending with $e_1e_2$, $f_1f_2$,
$f_1e_2$ and in the relation (9) terms ending with $e_1f_2$,
$f_1e_2$. Then we multiply both sides of (6) by
$$
(q^{H_1+H_2-1}+q^{-H_1-H_2+1})^{-1}
(q^{H_1+H_2+1}+q^{-H_1-H_2-1})^{-1}
$$
and both sides of (9) by
$$
(q^{H_1-H_2-1}+q^{-H_1+H_2+1})^{-1}
(q^{H_1-H_2+1}+q^{-H_1+H_2-1})^{-1} .
$$
After that we cancel in (6) terms ending with $e_1f_2$ and in (9)
terms ending with $e_1e_2$ and $f_1f_2$.
Now we have in (6) and (9) only the terms
ending with $e_1f_1$, $e_2f_2$ and terms without any $e_jf_k$, $j\ne k$.
We replace $e_1f_1$, $e_2f_2$ by the expressions following from
the expression for the Casimir elements, that is, by
$c_i-(q^{2H_i-1}+q^{-2H_i+1})/(q-q^{-1})^2$, respectively,
and multiply both sides of the both relations by
$$
(q^{H_1+H_2-1}+ q^{-H_1-H_2+1})
(q^{H_1+H_2+1}+ q^{-H_1-H_2-1})
(q^{H_1+H_2}+ q^{-H_1-H_2})
$$        $$
\times
(q^{H_1-H_2-1}+ q^{-H_1+H_2+1})
(q^{H_1-H_2+1}+ q^{-H_1+H_2-1})
(q^{H_1-H_2}+ q^{-H_1+H_2}).
$$
We obtain the relations both sides of which are
cancelled. Theorem is proved.
\medskip

\noindent
{\sf 5. PROPERTIES OF REPRESENTATIONS OF $U'_q({\rm so}_4)$}
\medskip

We assume everywhere below that $q$ is not a root of unity.

Our aim is to obtain irreducible finite dimensional representations of
$U'_q({\rm so}_4)$ by using the homomorphism of Theorem 2. But before we
need some statements on such representations of $U'_q({\rm so}_4)$.

Let $U'_q({\rm so}_3)_1$ and $U'_q({\rm so}_3)_2$ denote the subalgebras
of $U'_q({\rm so}_4)$ generated by $I_{21},I_{32}$ and by
$I_{32},I_{43}$, respectively. It is
known that the restriction of a finite dimensional representation $T$ of
$U'_q({\rm so}_4)$ to any of these subalgebras is a completely reducible
representation since any finite dimensional representation of
$U'_q({\rm so}_3)$ for $q$ not a root of unity is completely reducible
(see [29]). Let
$$
T\downarrow _{U'_q({\rm so}_3)_1} =R_1\oplus R_2\oplus \cdots \oplus R_k ,
\eqno (26)
$$    $$
T\downarrow _{U'_q({\rm so}_3)_2} =R'_1\oplus R'_2\oplus \cdots \oplus
R'_{k'} .   \eqno (27)
$$

Below we shall prove some assertions characterizing these decompositions.
Using the classification of irreducible finite dimensional representations
of $U'_q({\rm so}_3)$ for $q$ not a root of unity (see [22]), we can
state that each of irreducible representations $R_i$ and $R'_j$ in (26)
and (27) is a representation of the classical type or a representation
of the nonclassical type.
\medskip

{\bf Proposition 1.} {\it The decomposition (26), as well as the
decomposition (27), contains only irreducible representations of the
classical type or only irreducible representations of the nonclassical
type.}
\medskip

{\sl Proof.} Let us prove our proposition for the decomposition (26). It
follows from the results of
[30] that the operators $T(I_{41})$, $T(I_{42})$,
$T(I_{43})$ form a tensor operator transforming under the vector
(3-dimensional) representation of $U'_q({\rm so}_3)_1$. It is known from
[21] that tensor product of a classical (nonclassical) type
irreducible representation by the vector representation contains in the
decomposition classical (respectively nonclassical) type representations.
For this reason,
if, for example, the representation $R_1$ in (26) is of the classical type
and $| l_1,m_1\rangle$, $m_1=-l_1,-l_1+1,\cdots ,l_1$,
are basis elements of its representation subspace, then according to
Wigner--Eckart theorem (see [30]) for vector tensor operator $\{ T(I_{41}) ,
T(I_{42}) , T(I_{43})\}$ the vectors $T(I_{4i})|l_1,m_1\rangle$,
$i=1,\, 2,\, 3$, are linear combinations of vectors
belonging to subspaces of irreducible representations of
$U'_q({\rm so}_3)_1$ of the classical type. The same assertion is true
for vectors belonging to subspaces of representations of
the nonclassical type: that is, if $R_1$ is of the nonclassical type, then
the vectors $T(I_{4i})| l_1,m_1\rangle$ are linear combinations of
vectors belonging to subspaces of irreducible representations of
the nonclassical type.
Thus, if in the decomposition (26) there exists an
irreducible representation of the classical type, then acting upon vectors
of the corresponding subspace by the operators $T(I_{4i})$,
$i=1,\, 2,\, 3$, we obtain vectors belonging to subspaces on which
irreducible representations of the classical type
are realized. Since the representation $T$ of $U'_q({\rm so}_4)$ is
irreducible, then in this case the decomposition (26) contains only
irreducible representations of the classical type. If
the decomposition (26) does not contain an irreducible representation of the
classical type, then all representations in this decomposition are of the
nonclassical type. Proposition is proved.
\medskip

{\bf Proposition 2.} {\it Both decompositions (26) and (27) contain
irreducible representations
of the same type (classical or nonclassical).}
\medskip

{\sl Proof.} In order to prove this proposition we note (see [22])
that eigenvalues of the operator $R(I_{21})$ of an irreducible
representation $R$ of $U'_q({\rm so}_3)$ are of the form
${\rm i}[m]$, ${\rm i}[m+1],\cdots $ if $R$ is of the classical type, and of
the form $\pm [m]_+$, $\pm [m+1]_+, \cdots $ if $R$ is of the nonclassical
type, where
$$
[m]_+=\frac{q^m+q^{-m}}{q-q^{-1}}. \eqno (28)
$$

Let the decomposition (26) consist of irreducible representations of the
classical type.
Then eigenvalues of the operators $R_i(I_{21})$ are of the form
${\rm i}[m]$, ${\rm i}[m+1],\cdots $. We state that then the operators
$R_i(I_{32})$ can be
diagonalized and their eigenvalues are also of this form. Really, the
algebra
$U'_q({\rm so}_3)_1$ has the automorphism $\tau$ such that $\tau
(I_{21})=I_{32}$ and
$\tau (I_{32})=I_{21}$. Since ${\rm Tr}\, R_i(I_{21})={\rm Tr}\,
R_i(I_{32})=0$ (this equality characterizes (see [21])
irreducible representations of
the classical type), then the representations $R'_i=R_i\circ \tau$ are
of the classical type. Moreover, $R'_i\sim R_i$ (since up to
equivalence there exists a single irreducible representation of the
classical type with a fixed dimension).
Then the operators $R'_i(I_{21})=R_i(I_{32})$ are diagonalizable and the
spectrum of $R_i(I_{32})$ coincides with that of $R_i(I_{21})$. Our
statement concerning the operator $R_i(I_{32})$ is proved.

Now we consider the decomposition (27). The operator $T(I_{32})=\sum _i
\oplus R_i(I_{32})$ coincides with the operator
$T(I_{32})=\sum _i \oplus R'_i(I_{32})$. We have found a form of the
spectrum of the operator $T(I_{32})$. Now we can conclude that the
operators $R'_i(I_{32})$ have eigenvalues of the form
${\rm i}[m]$, ${\rm i}[m+1],\cdots $. This means that the irreducible
representations $R'_i$ in (27) are of the classical type.

If the decomposition (26) consists of irreducible representations of the
nonclassical type, the decomposition (27) consists of representations
of the same type. Really, if the
decomposition (27) would consist of representations of the classical type,
then conducting the above reasoning in the converse order
we would conclude that the decomposition (26) consists of representations
of the classical type. Proposition is proved.
\medskip

{\bf Corollary.} {\it If $T$ is an irreducible finite dimensional
representation of the algebra $U'_q({\rm so}_4)$, then both operators
$T(I_{21})$ and $T(I_{43})$ can be simulteneously
diagonalized and eigenvalues of both these operators are of the form
${\rm i}[m]$, ${\rm i}[m+1],\cdots $ or of the form $[m]_+$, $[m+1]_+,
\cdots $.}
\medskip

{\sl Proof.} The operators $T(I_{21})$ and $T(I_{43})$ can be simulteneously
diagonalized since the operators $R_i(I_{21})$ and $R'_j(I_{43})$ in the
decompositions (26) and (27) can
be simulteneously diagonalized. (Note that elements $I_{21}$ and $I_{43}$
are
commuting in $U'_q({\rm so}_4)$.) Since decompositions (26) and (27) consist
only of irreducible representations of the classical type or only
irreducible representations of the
nonclassical type, then the second assertion of the corollary follows.
Corollary is proved.
\medskip

\noindent
{\sf 6. IRREDUCIBLE REPRESENTATIONS OF $U'_q({\rm so}_4)$ OF
THE CLASSICAL TYPE}
\medskip

If $T$ is a representation of the algebra
${\hat U}_q({\rm sl}_2)^{\otimes 2,{\rm ext}}$
on a finite dimensional linear space ${\cal H}$, then the mapping
$$
R: U'_q({\rm so}_4)\to {\cal L}({\cal H})   \eqno (29)
$$
(where ${\cal L}({\cal H})$ is the space of linear operators on ${\cal H}$)
defined by the
composition $R=T\circ \phi$ (where $\phi$ is the homomorphism of Theorem 2)
is a representation of $U'_q({\rm so}_4)$. Let us consider the
representations
$$
R_{jj'}\equiv R^{(1,1)}_{jj'}=(T^{(1)}_j\otimes T^{(1)}_{j'})\circ \phi,
\ \ \ \ j,j'=0,\frac 12 , 1, \frac 32 ,\cdots ,   \eqno (30)
$$
of $U'_q({\rm so}_4)$, where $T^{(1)}_j\otimes T^{(1)}_{j'}$ are
irreducible representations of ${\hat U}_q({\rm sl}_2)^{\otimes 2,
{\rm ext}}$ from Theorem 1.

Using formulas for the representations
$T^{(1)}_j\otimes T^{(1)}_{j'}$ of
${\hat U}_q({\rm sl}_2)^{\otimes 2,{\rm ext}}$ from Theorem 1 and
the expressions (24) and (25) for $\phi (I_{i,i-1})$, $i=2,3,4$,
we find that
$$
R_{jj'}(I_{21})|k,l\rangle ={\rm i}[k+l]|k,l\rangle ,   \eqno (31)
$$    $$
R_{jj'}(I_{43})|k,l\rangle ={\rm i}[k-l]|k,l\rangle ,   \eqno (32)
$$    $$
R_{jj'}(I_{32})|k,l\rangle =\frac{1}{(q^{k+l}+q^{-k-l})(q^{k-l}+q^{-k+l})}
\times
$$      $$
\times \{ -(q^{j-l}+q^{-j+l}) [j'-l] |k,l+1\rangle +
(q^{j+l}+q^{-j-l}) [j'+l] |k,l-1\rangle +
$$     $$
+(q^{j'-k}+q^{-j'+k}) [j-k] |k+1,l\rangle -
(q^{j'+k}+q^{-j'-k}) [j+k] |k-1,l\rangle \} ,  \eqno (33)
$$
where numbers in square brackets are corresponding $q$-numbers and
$|k,l\rangle$ denote the
basis vetor
$$
|k,l\rangle \equiv |j,k\rangle \otimes |j',l\rangle
$$
of the space ${\cal H}_j\otimes {\cal H}_{j'}$ of the representation
$T^{(1)}_j\otimes T^{(1)}_{j'}$ of
${\hat U}_q({\rm sl}_2)^{\otimes 2,{\rm ext}}$.
\medskip

{\sl Remark:} Taking instead of $T^{(1)}_j\otimes T^{(1)}_{j'}$
the irreducible representations with other values of the operators
corresponding to the elements $x_1^k\otimes x_2^s$, $k,s=1,2,3$ (see
Theorem 1), we would obtain representations of $U'_q({\rm so}_4)$
equivalent to $R_{jj'}$.
\medskip

The representation $R_{jj'}$ of $U'_q({\rm so}_4)$ is equivalent to the
representation $T_{rs}$, $r=j+j'$, $s=j-j'$, from [24]
which in the $U'_q({\rm so}_3)$ basis
$$
|j'',m\rangle ,\ \ \ \ |s|\le j''\le r,\ \ \ \ m=-j'',-j''+1,\cdots ,j'',
$$
is given by the formulas
$$
T_{rs}(I_{21}) |j'',m\rangle ={\rm i}[m] |j'',m\rangle  ,   \eqno (34)
$$          $$
T_{rs}(I_{32}) |j'',m\rangle =\frac{1}{q^m+q^{-m}} ([j''-m]|j'',m+1\rangle
-[j''+m] |j'',m-1\rangle ), \eqno (35)
$$           $$
T_{rs}(I_{43})|j'',m\rangle ={\rm i}\frac{[r+1][s][m]}{[j''][j''+1]}
|j'',m\rangle  + \frac{[r-j''][j''+s+1]}{[j''+1][2j''+1]} |j''+1,m\rangle  -
$$      $$
- \frac{[r+j''+1][j''-s][j''-m][j''+m]}{[j''][2j''+1]} |j''-1,m\rangle
\eqno (36)
$$
(note that our basis elements $|j'',m\rangle $ differ from the basis
elements in formula (19) of [24] by the appropriate multipliers; for
this reason, formula (36) differs from formula (19) in [24]).
It was shown in [24] that under diagonalization of the operator
$T_{rs}(I_{43})$ we obtain a new basis
$\{ |x,m\rangle \}$ on which the operators $T_{rs}(I_{21})$,
$T_{rs}(I_{43})$, $T_{rs}(I_{32})$ are given by formulas
$$
T_{rs}(I_{21}) |x,m\rangle ={\rm i}[m]|x,m\rangle ,\ \ \ \
T_{rs}(I_{43}) |x,m\rangle ={\rm i}[x] |x,m\rangle \eqno (37)
$$
and by formula (33) in [24]. Under the renotation of the basis elements
and
multiplying them by the appropriate multipliers we obtain formulas for the
representation $R_{jj'}$ obtained above. This proves the equivalence stated
above.

The irreducible representations $R_{jj'}$ are called representations
of the classical type (see [14]). If $q\to 1$, then the operators
$R_{jj'}(I_{i,i-1})$, $i=2,3,4$, tend to the corresponding operators
of the irreducible representations of the Lie algebra ${\rm so}_4$.
\bigskip

\noindent
{\sf 7. IRREDUCIBLE REPRESENTATIONS OF $U'_q({\rm so}_4)$
OF THE NONCLASSICAL TYPE}
\medskip

Now we apply the method of the previous section to the irreducible
representations $T_j^{(-{\rm i})}\otimes T^{(1)}_{j'}$ of the algebra
${\hat U}_q({\rm sl}_2)^{\otimes 2,{\rm ext}}$ with $j$ half-integral and
$j'$ integral or with $j$ integral and $j'$ half-integral. Then we obtain
the representations
$$
R_{jj'}^{(-{\rm i},1)} = (T_j^{(-{\rm i})}\otimes T^{(1)}_{j'} )\circ \phi ,
\eqno (38)
$$
of $U'_q({\rm so}_4)$ with
$$
j=\frac 12 ,\frac 32 ,\frac 52 ,\cdots , \ \  j'=0,1,2,\cdots \ \ \ \
{\rm or} \ \ \ \  j=0,1,2,\cdots ,\ \  j'=\frac 12 ,\frac 32 ,\frac 52 ,
\cdots .  \eqno (39)
$$
Using formulas for the representations $T_j^{(-{\rm i})}$ and
$T^{(1)}_{j'}$ of $U_q({\rm sl}_2)$ and the expressions
(24) and (25) for $\phi (I_{i,i-1})$, $i=2,3,4$, we find that
$$
R_{jj'}^{(-{\rm i},1)} (I_{21})|k,l\rangle = [k+l]_+ |k,l\rangle ,
\eqno (40)
$$      $$
R_{jj'}^{(-{\rm i},1)} (I_{43}) |k,l\rangle =[k-l]_+ |k,l\rangle ,
\eqno (41)
$$       $$
R_{jj'}^{(-{\rm i},1)} (32) |k,l\rangle =\frac{1}{[k+l][k-l](q-q^{-1})}
\times
$$      $$
\times \left\{ -{\rm i} [j-l][j'-l] |k,l+1\rangle +
{\rm i}[j+l][j'+l] |k,l-1\rangle \right. -
$$     $$ \left.
- {\rm i} [j'-k][j-k] |k+1,l\rangle +
{\rm i}[j'+k][j+k] |k-1,l\rangle  \right\} ,    \eqno (42)
$$
where the numbers in square brackets are the corresponding $q$-numbers,
$[a]_+$ are defined in (38),
and $|k,l\rangle $ denote the basis vectors
$$
|k,l\rangle = |j,k\rangle \otimes |j',l\rangle     \eqno (43)
$$
of the space ${\cal H}_j\otimes {\cal H}_{j'}$ of the representation
$T_j^{(-{\rm i})}\otimes T^{(1)}_{j'}$ of
${\hat U}_q({\rm sl}_2)^{\otimes 2,{\rm ext}}$. Note that both $j+j'$ and
$k+l$ are half-integral.

In this case we have the equalities
$[k+l]_+=[-k-l]_+$, $[k-l]_+=[-k+l]_+$ and for this
reason the operators $R_{jj'}^{(-{\rm i},1)} (I_{21})$ and
$R_{jj'}^{(-{\rm i},1)} (I_{43})$ have multiple eigenvalues. Namely,
all pairs of the
vectors $|k,l\rangle $ and $|-k,-l\rangle $  have the same eigenvalues:
$$
R_{jj'}^{(-{\rm i},1)} (I_{21}) |k,l\rangle =[k+l]_+|k,l\rangle ,\ \ \ \
R_{jj'}^{(-{\rm i},1)} (I_{43}) |k,l\rangle =[k-l]_+|k,l\rangle ,
$$         $$
R_{jj'}^{(-{\rm i},1)} (I_{21}) |-k,-l\rangle =[k+l]_+ |-k,-l\rangle ,\ \ \
\
R_{jj'}^{(-{\rm i},1)} (I_{43}) |-k,-l\rangle =[k-l]_+ |-k,-l\rangle .
$$
The representation $R_{jj'}^{(-{\rm i},1)}$ is reducible.
In order to show this we distinguish two case:
\medskip

(a) $j$ is half-integral and $j'$ integral;

(b) $j$ is integral and $j'$ half-integral.
\medskip

\noindent
We consider first the case (a). In order to decompose
$R_{jj'}^{(-{\rm i},1)}$
into irreducible constituents we choose a new basis
in the representation space consisting of the vectors
$$
|k,l\rangle ^+ = |k,l\rangle +(-1)^l |-k,-l\rangle  ,\ \ \ \ \ k>0,
$$      $$
|k,l\rangle ^- =|k,l\rangle +(-1)^{l+1} |-k,-l\rangle ,\ \ \ \ \ k>0.
$$
Using formulas (40)--(42) we easily find that
$$
R_{jj'}^{(-{\rm i},1)} (I_{21}) |k,l\rangle ^{\pm}=[k+l]_+
|k,l\rangle ^{\pm} ,\ \ \
R_{jj'}^{(-{\rm i},1)} (I_{43}) |k,l\rangle ^{\pm}=[k-l]_+
|k,l\rangle ^{\pm},   \eqno (44)
$$         $$
R_{jj'}^{(-{\rm i},1)} (I_{32}) |k,l\rangle ^+= \frac{1}
{[k+l][k-l](q-q^{-1})} \times
$$        $$
\times \left\{ -{\rm i} [j-l][j'-l] |k,l+1\rangle ^++
{\rm i}[j+l][j'+l] |k,l-1\rangle ^+ \right.  -
$$     $$     \left.
-{\rm i} [j'-k][j-k] |k+1,l\rangle ^++
{\rm i} [j'+k][j+k] |k-1,l\rangle ^+ \right\} ,\ \ k\ne \frac 12 ,
 \eqno (45)
$$      $$
R_{jj'}^{(-{\rm i},1)} (I_{32}) | 1/2 ,l\rangle ^+ =\frac{-1}
{ [l+1/2] [l-1/2] (q-q^{-1})} \times
$$        $$
\times \left\{ -{\rm i} [j-l][j'-l] | 1/2 ,l+1\rangle ^+ +
{\rm i} [j+l] [j'+l] | 1/2 ,l-1\rangle ^+ \right.  -
$$     $$ \left.
-{\rm i} [j'-1/2][j-1/2] |3/2,l\rangle ^+ +
{\rm i} [j'+1/2] [j+1/2](-1)^l |1/2,-l\rangle ^+ \right\} .
\eqno (46)
$$
The operator $R_{jj'}^{(-{\rm i},1)} (I_{32})$ acts on the vectors
$|k,l\rangle ^-$ by formulas (45) and (46) if to replace all
$|k,l\rangle ^+$ by the corresponding $|k,l\rangle ^-$ and in (46)
to replace $(-1)^l$ by $(-1)^{l+1}$.

Let ${\cal H}_1$ and ${\cal H}_2$ be the subspaces of the representation
space of $R_{jj'}^{(-{\rm i},1)}$ spanned by the vectors $|k,l\rangle ^+$
and by the vectors $|k,l\rangle ^-$, respectively. We see from the above
formulas that these subspaces are invariant with respect to the
representation $R_{jj'}^{(-{\rm i},1)}$. We denote the restrictions of
$R_{jj'}^{(-{\rm i},1)}$ to ${\cal H}_1$ and ${\cal H}_2$ by
$R_{jj'}^{(+,+,+)}$ and $R_{jj'}^{(+,+,-)}$, respectively. Note that
$$
R_{jj'}^{(+,+,+)} (I_{32}) | 1/2 ,0\rangle ^+ =\frac{1}{[1/2]^{2}
(q-q^{-1})} \left\{ -{\rm i} [j][j'] | 1/2 ,1\rangle ^+ +
{\rm i} [j] [j'] | 1/2 ,-1\rangle ^+ \right. -
$$     $$     \left.
-{\rm i} [j'-1/2][j-1/2] |3/2,0\rangle ^+ +
{\rm i} [j'+1/2] [j+1/2] |1/2,0\rangle ^+ \right\} ,
\eqno (47)
$$
that is, the operator $R_{jj'}^{(+,+,+)} (I_{32})$ has nonzero diagonal
element
$$
{}^+\langle 1/2 ,0\, | R_{jj'}^{(+,+,+)} (I_{32}) \, | 1/2 ,
0\rangle ^+ = {\rm i } \frac{[j'+1/2][j+1/2]}{[1/2]^2(q-q^{-1})} .
$$
In the same way it is shown that
$$
{}^-\langle 1/2 ,0\, | R_{jj'}^{(+,+,-)} (I_{32}) \, | 1/2 ,
0\rangle ^- = -{\rm i } \frac{[j'+1/2][j+1/2]}{[1/2]^2(q-q^{-1})} .
$$

Now we consider the case (b). The new basis of the representation
space is
$$
|k,l\rangle ^+ = |k,l\rangle +(-k)^l |-k,-l\rangle  ,\ \ \ \ \ l>0,
$$      $$
|k,l\rangle ^- =|k,l\rangle +(-1)^{k+1} |-k,-l\rangle ,\ \ \ \ \ l>0.
$$
Then formulas for $R_{jj'}^{(-{\rm i},1)} (I_{21})$ and
$R_{jj'}^{(-{\rm i},1)} (I_{43})$ are the same as in case (a).
Formula (45) is the same, but we have consider it for $l\ne \frac 12$.
Instead of (46) we have the formula
$$
R_{jj'}^{(-{\rm i},1)} (I_{32}) |k, 1/2 \rangle ^+ =\frac{1}
{ [k+1/2] [k-1/2] (q-q^{-1})} \times
$$        $$
\times \left\{ -{\rm i} [j-1/2][j'-1/2] |k, 3/2 \rangle ^+ +
{\rm i} [j+1/2] [j'+1/2] (-1)^k|-k, 1/2 \rangle ^+ \right. -
$$     $$ \left.
-{\rm i} [j'-k][j-k] |k+1,1/2\rangle ^+ +
{\rm i} [j'+k] [j+k] |k-1,1/2\rangle ^+ \right\} .
\eqno (48)
$$
The operator $R_{jj'}^{(-{\rm i},1)}(I_{32})$ acts on the vectors
$|k,l\rangle ^-$ by formulas (45) and (48) if to replace all
$|k,l\rangle ^+$ by $|k,l\rangle ^-$, respectively, and in (48)
to replace $(-1)^k$ by $(-1)^{k+1}$.

We denote the subspaces spanned by the vectors $|k,l\rangle ^+$
and by the vectors $|k,l\rangle ^-$  as ${\cal H}_1$ and
${\cal H}_2$, respectively. As we see from above formulas,
these subspaces are invariant with respect to the
representation $R_{jj'}^{(-{\rm i},1)}$. We denote the
corresponding subrepresentations by
$R_{jj'}^{(+,+,-)}$ and $R_{jj'}^{(+,+,+)}$, respectively. In
particular, we have
$$
R_{jj'}^{(+,+,-)} (I_{32}) |0, 1/2 \rangle ^+ =\frac{-1}{[1/2]^{2}
(q-q^{-1})} \left\{ -{\rm i} [j-1/2][j'-1/2] |0, 3/2 \rangle ^+ +
\right.
$$    $$
+{\rm i} [j+1/2] [j'+1/2] |0, 1/2 \rangle ^+ -
\left.
{\rm i} [j'][j] |1,1/2\rangle ^+ +
{\rm i} [j'] [j]|-1,1/2\rangle ^+ \right\} ,
\eqno (49)
$$
that is, the operator $R_{jj'}^{(+,+,-)} (I_{32})$ has nonzero
diagonal element
$$
{}^+\langle 0, 1/2 \, | R_{jj'}^{(+,+,-)} (I_{32}) \, | 0, 1/2
\rangle ^+ = -{\rm i } \frac{[j'+1/2][j+1/2]}{[1/2]^2(q-q^{-1})} .
$$
For the operator $R_{jj'}^{(+,+,+)}(I_{32})$ we have
$$
{}^-\langle 0, 1/2 \, | R_{jj'}^{(+,+,+)} (I_{32}) \, |0, 1/2
\rangle ^- = {\rm i } \frac{[j'+1/2][j+1/2]}{[1/2]^2(q-q^{-1})} .
$$

Now we consider in the same way the representations
$$
R_{jj'}^{(-{\rm i},-1)}= (T_j^{(-{\rm i})}\otimes T^{(-1)}_{j'})
\circ \phi
$$
of $U'_q({\rm so}_4)$. As a result, we obtain for these representations
formulas (40)--(42) in which right hand sides of (40) and (41) are
multiplied by $-1$ and the right hand side of (42) is left
without any change. These representations are reducible and we have
$$
R_{jj'}^{(-{\rm i},-1)} = R_{jj'}^{(-,-,+)}\oplus R_{jj'}^{(-,-,-)} ,
$$
where the representations $R_{jj'}^{(-,-,+)}$ and $R_{jj'}^{(-,-,-)}$
act on such subspaces as in the previous case, with those difference
that
$$
R_{jj'}^{(-,-,\pm )}(I_{21})=- R_{jj'}^{(+,+,\pm )}(I_{21}),\ \ \ \ \
R_{jj'}^{(-,-,\pm )} (I_{43})=- R_{jj'}^{(+,+,\pm )} (I_{43}) ,
$$
that is,
$$
R_{jj'}^{(-,-,\pm )}(I_{21})|k,l\rangle ^\pm =-[k+l]_+
|k,l\rangle ^\pm ,\ \ \ \
R_{jj'}^{(-,-,\pm )}(I_{43})|k,l\rangle ^\pm =-[k-l]_+
|k,l\rangle ^\pm  .
$$
The operators $R_{jj'}^{(-,-,\pm )}(I_{32})$ coincide with the
corresponding operators $R_{jj'}^{(+,+,\pm )}(I_{32})$.

Similarly, the representation
$R_{jj'}^{({\rm i},1)}= (T_j^{({\rm i})}\otimes T^{(1)}_{j'})
\circ \phi $ of $U'_q({\rm so}_4)$ is reducible and
decomposes as
$$
R_{jj'}^{({\rm i},1)} = R_{jj'}^{(-,-,+)}\oplus R_{jj'}^{(-,-,-)} ,
$$
where $R_{jj'}^{(-,-,\pm )}$ are such as above, and the representation
$R_{jj'}^{({\rm i},-1)}= (T_j^{({\rm i})}\otimes T^{(-1)}_{j'})
\circ \phi $ of $U'_q({\rm so}_4)$ is reducible and
decomposes as
$$
R_{jj'}^{({\rm i},-1)} = R_{jj'}^{(+,+,+)}\oplus R_{jj'}^{(+,+,-)} .
$$

Let us consider the representation
$$
R_{j'j}^{(1,{\rm i})}= (T^{(1)}_{j'} \otimes T_{j}^{({\rm i})})
\circ \phi
$$
of $U'_q({\rm so}_4)$ (for convenience we take first the index $j'$
and then $j$). We have
$$
R_{j'j}^{(1,{\rm i})} (I_{21})|l,k\rangle = -[k+l]_+ |l,k\rangle ,\ \ \ \ \
R_{j'j}^{(1,{\rm i})} (I_{43}) |l,k\rangle =[k-l]_+ |l,k\rangle ,
$$       $$
R_{j'j}^{(1,{\rm i})} (I_{32}) |l,k\rangle =
\frac{1}{[k+l][l-k](q-q^{-1})} \times
$$      $$
\times \left\{ -{\rm i} [j-l][j'-l] |l+1,k\rangle +
{\rm i}[j+l][j'+l] |l-1,k\rangle  \right. -
$$     $$   \left.
-{\rm i} [j'-k][j-k] |l,k+1\rangle +
{\rm i}[j'+k][j+k] |l,k-1\rangle  \right\}   .
$$
In order to have a similarity with formulas (40)--(42) we denote
the vectors $|l,k\rangle$ by $|k,l\rangle$ and the representation
$R_{j'j}^{(1,{\rm i})}$ by ${\hat R}_{jj'}^{({\rm i},1)}$.

The operators ${\hat R}_{jj'}^{({\rm i},1)} (I_{21})$ and
${\hat R}_{jj'}^{({\rm i},1)} (I_{43})$ have multiple common eigenvectors
$|k,l\rangle$ and $| -k,-l\rangle$. The representation
${\hat R}_{jj'}^{({\rm i},1)} $ is reducible and decomposes into
irreducible constituents in the same way as in the previous case.
The representation ${\hat R}_{jj'}^{({\rm i},1)}$ is
analysed as the representation
$R_{jj'}^{({-\rm i},1)}$ above and we have the following result.
The representation ${\hat R}_{jj'}^{({\rm i},1)}$ is reducible
and decomposes into irreducible subrepresentations as
$$
{\hat R}_{jj'}^{({\rm i},1)}=R^{(-,+,+)}_{jj'}\oplus
R^{(-,+,-)}_{jj'} ,
$$
where the representations $R^{(-,+,\pm )}_{jj'}$
differ from the representations $R^{(+,+,\pm )}_{jj'}$,
respectively, only by the operator $R^{(-,+,\pm )}_{jj'}(I_{21})$
and
$$
R^{(-,+,\pm )}_{jj'}(I_{21})=-R^{(+,+,\pm )}_{jj'}(I_{21}).
$$

Similarly, the representation
$$
R_{j'j}^{(1,-{\rm i})}\equiv  {\hat R}_{jj'}^{(-{\rm i},1)}
=(T^{(1)}_{j'} \otimes T_{j}^{(-{\rm i})})
\circ \phi
$$
of $U'_q({\rm so}_4)$ is reducible and decomposes into irreducible
components as
$$
{\hat R}_{jj'}^{(-{\rm i},1)}=R^{(+,-,+)}_{jj'}\oplus
R^{(+,-,-)}_{jj'} ,
$$
where the representations $R^{(+,-,\pm )}_{jj'}$
differ from the representations $R^{(+,+,\pm )}_{jj'}$,
respectively, only by the operator $R^{(+,-,\pm )}_{jj'}(I_{43})$
and
$$
R^{(+,-,\pm )}_{jj'}(I_{43})=-R^{(+,+,\pm )}_{jj'}(I_{43}) .
$$

We do not consider other representations
$R_{j'j}^{(\pm 1,\pm {\rm i})}
=(T^{(\pm 1)}_{j'} \otimes T_{j}^{(\pm {\rm i})})$ since
they do not give new irreducible representations of
$U'_q({\rm so}_4)$.

Thus, for every values of $j$ and $j'$ such that
$$
j=0,1,2,\cdots ,\ \ j'={\frac 12},{\frac 32}, {\frac 53},\cdots
\ \ \ \ \ {\rm or}\ \ \ \ \
j={\frac 12},{\frac 32}, {\frac 53},\cdots ,\ \ j'=0,1,2,\cdots
$$
we constructed 8 representations
$R^{(\varepsilon _1,\varepsilon _2,\varepsilon _3)}_{jj'}$,
$\varepsilon _i=\pm$. These representations
act on the linear space ${\cal H}$ with the basis
$$
|k,l\rangle ,\ \ \ \ k=j,j-1,\cdots ,{\frac 12},\ \ \
l=j',j'-1,\cdots ,-j',
$$
if $j'$ is integral and with the basis
$$
|k,l\rangle ,\ \ \ \ k=j,j-1,\cdots ,-j,\ \ \
l=j',j'-1,\cdots ,{\frac 12},
$$
if $j$ is integral. The representations are given by the formulas
$$
R^{(\varepsilon _1,\varepsilon _2,\varepsilon _3)}_{jj'}(I_{21})
|k,l\rangle = \varepsilon _1[k+l]_+|k,l\rangle , \eqno (50)
$$      $$
R^{(\varepsilon _1,\varepsilon _2,\varepsilon _3)}_{jj'}(I_{43})
|k,l\rangle = \varepsilon _2 [k-l]_+ |k,l\rangle , \eqno (51)
$$       $$
R^{(\varepsilon _1,\varepsilon _2,\varepsilon _3)}_{jj'} (I_{32})
|k,l\rangle  =\frac{1}{[k+l][k-l](q-q^{-1})} \{
-{\rm i}[j'-l][j-l]|k,l+1\rangle +
$$       $$
+ {\rm i}[j'+l][j+l] |k,l-1\rangle -{\rm i} [j'-k][j-k]
|k+1,l\rangle + {\rm i} [j'+k][j+k] |k-1,l\rangle \},
\eqno (52)
$$
where $k\ne {\frac 12}$ if $j$ is half-integral and
$l\ne {\frac 12}$ if $j'$ is half-integral, and by
$$
R^{(\varepsilon _1,\varepsilon _2,\varepsilon _3)}_{jj'}(I_{32})
| 1/2 ,l\rangle  =\frac{1}
{ [l+1/2] [l-1/2] (q-q^{-1})}
\{ -{\rm i} [j-l][j'-l] |1/2 ,l+1\rangle  +
$$    $$
+{\rm i} [j+l] [j'+l] |\frac 12 ,l-1\rangle  -
{\rm i} [j'-1/2][j-1/2] |3/2,l\rangle  +
{\rm i} [j'+1/2] [j+1/2]\varepsilon _3 (-1)^l |1/2,-l\rangle \}
\eqno (53)
$$
if $j$ is half-integral and by
$$
R^{(\varepsilon _1,\varepsilon _2,\varepsilon _3)}_{jj'}(I_{32})
|k,1/2 \rangle  =\frac{1}
{ [k+1/2] [k-1/2] (q-q^{-1})}
\{ -{\rm i} [j-1/2][j'-1/2] |k, 3/2 \rangle  +
$$      $$
+{\rm i} [j+1/2] [j'+1/2]\varepsilon _3 (-1)^k|-k, 1/2 \rangle  -
{\rm i} [j'-k][j-k] |k+1,1/2\rangle  +
{\rm i} [j'+k] [j+k] |k-1,1/2\rangle \}
\eqno (54)
$$
if $j'$ is half-integral.

It is seen from formulas for representations that the representations
$R^{(\varepsilon _1,\varepsilon _2,\varepsilon _3)}_{jj'}$ and
$R^{(\varepsilon _1,-\varepsilon _2,\varepsilon _3)}_{j'j}$
are equivalent. The equivalence operator $A$ is given by the formula
$A|k,l\rangle =|l,k\rangle$. For this reason, we consider the
representations
$R^{(\varepsilon _1,\varepsilon _2,\varepsilon _3)}_{jj'}$ only
for $j\ge j'$.
\medskip

{\bf Theorem 3.} {\it The representations}
$R^{(\varepsilon _1,\varepsilon _2,\varepsilon _3)}_{jj'}$, $j\ge j'$,
{\it are irreducible and pairwise nonequivalent.}
\medskip

{\sl Proof.} Irreducibility of these representations
will be proved in section 8. In order to prove their
pairwise nonequivalence we note that finite dimensional
representations $T$ and $T'$ of $U'_q({\rm so}_4)$ cannot
be equivalent if at least for one pair $T(I_{i,i-1})$,
$T'(I_{i,i-1})$, $i=2,3,4$, we have
${\rm Tr}\, T(I_{i,i-1})\ne {\rm Tr}\, T'(I_{i,i-1})$.
For this reason, the representations
$R^{(\varepsilon _1,\varepsilon _2,\varepsilon _3)}_{jj'}$ and
$R^{(\varepsilon '_1,\varepsilon '_2,\varepsilon '_3)}_{ss'}$
with $(\varepsilon _1,\varepsilon _2,\varepsilon _3)\ne
(\varepsilon '_1,\varepsilon '_2,\varepsilon '_3)$ or/and with
$(j,j')\ne (s,s')$ are not
equivalent.
\medskip

Irreducible representations
$R^{(\varepsilon _1,\varepsilon _2,\varepsilon _3)}_{jj'}$
are called {\it representations of the nonclassical type}. They have no
classical analogue. Their main property is that the operators
$R^{(\varepsilon _1,\varepsilon _2,\varepsilon _3)}_{jj'}(I_{i,i-1})$,
$i=2,3,4$, have a nonzero trace. Note that there exist 8 nontrivial
one-dimensional representations of the nonclassical type. They
coincide with the representations
$R^{(\varepsilon _1,\varepsilon _2,\varepsilon _3)}_{{\frac 12} 0}$.
\bigskip

\noindent
{\sf 8. THE CLASSIFICATION THEOREM}
\medskip

The main aim of this section is to prove that constructed
above irreducible representations of the classical type
and of the nonclassical type exhaust all irreducible
finite dimensional representations of the algebra
$U'_q({\rm so}_4)$. We also prove that pairs of the
irreducible representations of Theorem 3 are not
equivalent. But first we study some auxiliary operators.

Let $R$ be a finite dimensional representation of
$U'_q({\rm so}_4)$ on a linear vector space ${\cal H}$. We
suppose that the operators $R(I_{21})$ and $R(I_{43})$
have eigenvalues only of the classical type, that is, of
the form ${\rm i} [m]$, where $[m]$ means a $q$-number.
Let $|k,l\rangle$ be an eigenvector such that
$$
R(I_{21}) |k,l\rangle ={\rm i}[k+l]|k,l\rangle ,\ \ \ \
R(I_{43}) |k,l\rangle ={\rm i}[k-l]|k,l\rangle .
$$
We associate with this eigenvector the operators
$$
X_1^{(k,l)}=-R(I_{41}) + q^{-2k}R(I_{32})-{\rm i}
q^{-k-l+1/2}R(I_{42}) -{\rm i}q^{-k+l-1/2}R(I_{31}) ,
\eqno (55)
$$     $$
X_2^{(k,l)}=-R(I_{41}) + q^{2k}R(I_{32})+{\rm i}q^{k+l+1/2}R(I_{42})
+{\rm i} q^{k-l-1/2}R(I_{31})  ,
\eqno (56)
$$     $$
X_3^{(k,l)}=R(I_{41}) + q^{-2l}R(I_{32})+{\rm i}
q^{-k-l+1/2}R(I_{42}) -{\rm i}q^{k-l-1/2}R(I_{31}) ,
\eqno (57)
$$     $$
X_4^{(k,l)}=R(I_{41}) + q^{2l}R(I_{32})-{\rm i}
q^{k+l+1/2}R(I_{42}) +{\rm i}q^{-k+l-1/2}R(I_{31}) .
\eqno (58)
$$

Note that below we shall use explicit form of the
operators $X_1^{(k,l)}$, $X_2^{(k,l)}$,
$X_3^{(k,l)}$, $X_4^{(k,l)}$ under action on the
vector $|k,l\rangle$. In this case we may consider
formulas (55)--(58) as a system of linear equations
with unknown vectors $R(I_{32})|k,l\rangle$,
$R(I_{31})|k,l\rangle$, $R(I_{42})|k,l\rangle$,
$R(I_{41})|k,l\rangle$. The determinant of this
system can be easily calculated:
$$
{\rm det} \left(
\matrix{
-1 & q^{-2k} & -{\rm i}q^{-k-l+1/2} &
          -{\rm i}q^{-k+l-1/2} \cr
-1 & q^{2k} & {\rm i}q^{k+l+1/2} &
          {\rm i}q^{k-l-1/2} \cr
1 & q^{-2l} & {\rm i}q^{-k-l+1/2} &
          -{\rm i}q^{k-l-1/2} \cr
1 & q^{2l} & -{\rm i}q^{k+l+1/2} &
                  {\rm i}q^{-k+l-1/2} } \right)
=(q^{k+l}+q^{-k-l})(q^{k-l}+q^{-k+l}).
$$
If $q$ is not a root of unity, then this
determinant does not vanish for any integral or
half-integral $k$ and $l$. This means that
the system of above equations can be solved
and we can find how the operators
$R(I_{32})$, $R(I_{31})$, $R(I_{42})$ and
$R(I_{41})$ act on the vector $|k,l\rangle$.
This reasoning will be used below.
\medskip

{\bf Lemma 1.} {\it The vectors} $X_i^{(k,l)}|k,l\rangle$
{\it are eigenvectors of the operators
$R(I_{21})$ and $R(I_{43})$:}
$$
R(I_{21}) (X_1^{(k,l)}|k,l\rangle )= {\rm i}
[k+l+1] (X_1^{(k,l)}|k,l\rangle ), \ \ \
R(I_{43}) (X_1^{(k,l)}|k,l\rangle )= {\rm i}
[k-l-1] (X_1^{(k,l)}|k,l\rangle ),
$$         $$
R(I_{21}) (X_2^{(k,l)}|k,l\rangle )= {\rm i}
[k+l-1] (X_2^{(k,l)}|k,l\rangle ), \ \ \
R(I_{43}) (X_2^{(k,l)}|k,l\rangle )= {\rm i}
[k-l+1] (X_2^{(k,l)}|k,l\rangle ),
$$       $$
R(I_{21}) (X_3^{(k,l)}|k,l\rangle )= {\rm i}
[k+l+1] (X_3^{(k,l)}|k,l\rangle ), \ \ \
R(I_{43}) (X_3^{(k,l)}|k,l\rangle )= {\rm i}
[k-l+1] (X_3^{(k,l)}|k,l\rangle ),
$$         $$
R(I_{21}) (X_4^{(k,l)}|k,l\rangle )= {\rm i}
[k+l-1] (X_4^{(k,l)}|k,l\rangle ), \ \ \
R(I_{43}) (X_4^{(k,l)}|k,l\rangle )= {\rm i}
[k-l-1] (X_4^{(k,l)}|k,l\rangle ).
$$
{\sl Proof.} The lemma is proved by direct calculation using
the defining relations for the elements $I_{ij}\in
U'_q({\rm so}_4)$, $i>j$. For example, by using
relations (11)--(14) we have
$$
R(I_{21}) (X_1^{(k,l)}|k,l\rangle ) =
R(-qI_{41}+q^{-2k-1}I_{32}-{\rm i}q^{-k-l-1/2}I_{42} -
{\rm i}q^{k+l+1/2}I_{31})R(I_{21}) |k,l\rangle +
$$       $$
+R({\rm i}q^{-k-l}I_{41}-{\rm i}q^{l-k}I_{32}
+q^{1/2}I_{42}+q^{-2k-1/2}I_{31}) |k,l\rangle =
{\rm i}[k+l+1](X_1^{(k,l)}|k,l\rangle ) .
$$
Lemma is proved.

 Lemma 1 means that the operators $X_1^{(k,l)}$
and $X_3^{(k,l)}$ increase $k$ and $l$ in the eigenvalues of
$R(I_{21})$ and $R(I_{43})$, respectively, and
the operators $X_2^{(k,l)}$
and $X_4^{(k,l)}$ decrease these numbers in these eigenvalues.
Symbolically we write this in the form
$$
X_1:\ \ l\to l+1,\ \ \ \ \ \ X_2:\ \ l\to l-1,
$$     $$
X_3:\ \ k\to k+1,\ \ \ \ \ \ X_4:\ \ k\to k-1 .
$$
{\bf Lemma 2.} {\it The operators (55)--(58) have the
properties}
$$
X_3^{(k,l+1)} X_1^{(k,l)} |k,l\rangle =
X_1^{(k+1,l)} X_3^{(k,l)} |k,l\rangle ,
\ \ \ \ \
X_4^{(k,l-1)} X_2^{(k,l)} |k,l\rangle =
X_2^{(k-1,l)} X_4^{(k,l)} |k,l\rangle ,
$$     $$
X_4^{(k,l+1)} X_1^{(k,l)} |k,l\rangle =
X_1^{(k-1,l)} X_4^{(k,l)} |k,l\rangle ,
\ \ \ \ \
X_3^{(k,l-1)} X_2^{(k,l)} |k,l\rangle =
X_2^{(k+1,l)} X_3^{(k,l)} |k,l\rangle .
$$
{\sl Proof.} The first of these relations is proved
as follows. Using the expressions (55) and (57) for
$X_1$ and $X_3$ we express (by using relations
(11)--(15)) the elements
$X_3^{(k,l+1)} X_1^{(k,l)}$ and
$X_1^{(k+1,l)} X_3^{(k,l)}$ as a linear combination
of the basis elements from Poincar{\'e}--Birkhoff--Witt
theorem. As a result, we recieve the first relation.
Other relations are proved in the same way. Lemma
is proved.
\medskip

Lemma 2 means that the pairs of operators $X_1$ and
$X_3$, $X_2$ and $X_4$, $X_1$ and $X_4$, $X_2$
and $X_3$ (with appropriate upper indices)
commute under action on the vector $|k,l\rangle$.
\medskip

{\bf Lemma 3.} {\it The operators (55)--(58) have the
properties}
$$
X_2^{(k,l+1)} X_1^{(k,l)} |k,l\rangle =
(C'_4 -(q^{2l+1}+q^{-2l-1})C_4+[2l][2(l+1)])
|k,l\rangle ,
$$       $$
X_1^{(k,l-1)} X_2^{(k,l)} |k,l\rangle =
(C'_4 -(q^{2l-1}+q^{-2l+1})C_4+[2l][2(l-1)])
|k,l\rangle ,
$$       $$
X_4^{(k+1,l)} X_3^{(k,l)} |k,l\rangle =
(C'_4 +(q^{2k+1}+q^{-2k-1})C_4+[2k][2(k+1)])
|k,l\rangle ,
$$       $$
X_3^{(k-1,l)} X_4^{(k,l)} |k,l\rangle =
(C'_4 +(q^{2k-1}+q^{-2k+1})C_4+[2k][2(k-1)])
|k,l\rangle ,
$$
{\it where $C_4$ and $C'_4$ are the Casimir elements of
$U'_q({\rm so}_4)$ from section 2.}
\medskip

{\sl Proof} is given in the same way as that of Lemma 2,
taking into accout expressions for the Casimir elements.
\medskip

Lemma 3 can be used for evaluation of eigenvalues of Casimir elements
$C_4$ and $C'_4$ on irreducible representations $R$ when the operators
$R(I_{21})$ and $R(I_{43})$ have eigenvalues of the classical type.
An eigenvectors $| k,l\rangle$ of the operators $R(I_{21})$ and $R(I_{43})$
are called {\it weight vectors} of the representation $R$. A weight
vector $| j,j'\rangle$ is called {\it a highest weight vector} if
$$
X^{(j,j')}_1 |j,j'\rangle =0,\ \ \ \ \  X^{(j,j')}_3 | j,j'\rangle =0.
$$
If $R$ is an irreducible representation with classical type eigenvalues
of the operators $R(I_{21})$ and $R(I_{43})$, then we apply both sides
of the first and the third relations of Lemma 3 to
the vector of highest weight $|j,j'\rangle$.
The left hand sides send this vector to zero, and
then the right hand sides (equating to zero)
gives the following eigenvalues for $C_4$ and $C'_4$:
$$
C_4=[j+j'+1][j'-j]I,
$$    $$
C'_4=\{ (q^{2j+1}+q^{-2j-1})[j-j'][j+j'+1]-[2j][2j+2]\} I ,
$$
where $I$ is the unit operator in the representation space.
In particular, such eigenvalues have Casimir operators of the classical
type representation $R_{jj'}$.

Now let $R'$ be a finite dimensional representation of
$U'_q({\rm so}_4)$ on a linear space ${\cal H}'$. Suppose
that the operators $R'(I_{21})$ and $R'(I_{43})$
have eigenvalues only of the nonclassical type, that is, of
the form $\pm [m]_+$, where $[m]_+=(q^m+q^{-m})/(q-q^{-1})$
and $m$ are half-integral.
If $|k,l\rangle$ is an eigenvector such that
$$
R'(I_{21}) |k,l\rangle =\varepsilon _1[k+l]_+|k,l\rangle ,
\ \ \ \
R'(I_{43}) |k,l\rangle =\varepsilon _2[k-l]_+|k,l\rangle .
$$
then we associate with it the operators
$$
X_1^{(k,l)}=R'(I_{41}) + q^{-2k}R'(I_{32})-
q^{-k-l+1/2}R'(I_{42}) -q^{-k+l-1/2}R'(I_{31}) ,
\eqno (59)
$$     $$
X_2^{(k,l)}=R'(I_{41}) + q^{2k}R'(I_{32})- q^{k+l+1/2}R'(I_{42})
-q^{k-l-1/2}R'(I_{31})  ,
\eqno (60)
$$     $$
X_3^{(k,l)}=-R'(I_{41}) - q^{-2l}R'(I_{32})+
q^{-k-l+1/2}R'(I_{42}) + q^{k-l-1/2}R'(I_{31}) ,
\eqno (61)
$$     $$
X_4^{(k,l)}=-R'(I_{41}) - q^{2l}R'(I_{32}) +
q^{k+l+1/2}R'(I_{42}) +q^{-k+l-1/2}R'(I_{31}) .
\eqno (62)
$$

Below we shall consider relations (59)--(62)
as a system of linear equations. Determinant
of the matrix of this system is equal to
$$
{\rm det} \left(
\matrix{
1 & q^{-2k} & -q^{-k-l+1/2} &
          -q^{-k+l-1/2} \cr
1 & q^{2k} &  -q^{k+l+1/2} &
          -q^{k-l-1/2} \cr
-1 & -q^{-2l} & q^{-k-l+1/2} &
           q^{k-l-1/2} \cr
-1 & -q^{2l} & q^{k+l+1/2} &
                   q^{-k+l-1/2} } \right)
=(q^{k+l}-q^{-k-l})(q^{k-l}-q^{-k+l}).
$$
If $q$ is not a root of unity, then this
determinant does not vanish for any
half-integral $k\pm l$. Hence,
the system of above equations can be solved
and we can find how the operators
$R'(I_{32})$, $R'(I_{31})$, $R'(I_{42})$ and
$R'(I_{41})$ act on the vector $|k,l\rangle$.

Below we formulate three lemmas for these operators
analogous to Lemma 1--3. Proofs of these lemmas are
the same as in the case of Lemmas 1--3 and we omit
them.
\medskip

{\bf Lemma 4.} {\it The vectors} $X_i^{(k,l)}|k,l\rangle$
{\it are eigenvectors of the operators
$R'(I_{21})$ and $R'(I_{43})$:}
$$
R'(I_{21}) (X_1^{(k,l)}|k,l\rangle ){=}
\varepsilon _1 [k+l+1]_+ (X_1^{(k,l)}|k,l\rangle ), \
R'(I_{43}) (X_1^{(k,l)}|k,l\rangle ){=}
\varepsilon _2[k-l-1]_+ (X_1^{(k,l)}|k,l\rangle ),
$$         $$
R'(I_{21}) (X_2^{(k,l)}|k,l\rangle ){=}
\varepsilon _1 [k+l-1]_+ (X_2^{(k,l)}|k,l\rangle ),  \
R'(I_{43}) (X_2^{(k,l)}|k,l\rangle ){=}
\varepsilon _2 [k-l+1]_+ (X_2^{(k,l)}|k,l\rangle ),
$$       $$
R'(I_{21}) (X_3^{(k,l)}|k,l\rangle ){=}
\varepsilon _1 [k+l+1]_+ (X_3^{(k,l)}|k,l\rangle ), \
R'(I_{43}) (X_3^{(k,l)}|k,l\rangle ){=}
\varepsilon _2 [k-l+1]_+ (X_3^{(k,l)}|k,l\rangle ),
$$         $$
R'(I_{21}) (X_4^{(k,l)}|k,l\rangle ){=}
\varepsilon _1 [k+l-1]_+ (X_4^{(k,l)}|k,l\rangle ),  \
R'(I_{43}) (X_4^{(k,l)}|k,l\rangle ){=}
\varepsilon _2 [k-l-1]_+ (X_4^{(k,l)}|k,l\rangle ).
$$

{\bf Lemma 5.} {\it The operators (59)--(62) have the
properties}
$$
X_3^{(k,l+1)} X_1^{(k,l)} |k,l\rangle =
X_1^{(k+1,l)} X_3^{(k,l)} |k,l\rangle ,
\ \ \ \
X_4^{(k,l-1)} X_2^{(k,l)} |k,l\rangle =
X_2^{(k-1,l)} X_4^{(k,l)} |k,l\rangle ,
$$     $$
X_4^{(k,l+1)} X_1^{(k,l)} |k,l\rangle =
X_1^{(k-1,l)} X_4^{(k,l)} |k,l\rangle ,
\ \ \ \
X_3^{(k,l-1)} X_2^{(k,l)} |k,l\rangle =
X_2^{(k+1,l)} X_3^{(k,l)} |k,l\rangle .
$$

{\bf Lemma 6.} {\it For the operators (59)--(62) we have}
$$
X_2^{(k,l+1)} X_1^{(k,l)} |k,l\rangle =
(C'_4 -(q^{2l+1}+q^{-2l-1})C_4+[2l][2(l+1)])
|k,l\rangle ,
$$       $$
X_1^{(k,l-1)} X_2^{(k,l)} |k,l\rangle =
(C'_4 -(q^{2l-1}+q^{-2l+1})C_4+[2l][2(l-1)])
|k,l\rangle ,
$$       $$
X_4^{(k+1,l)} X_3^{(k,l)} |k,l\rangle =
(C'_4 -(q^{2k+1}+q^{-2k-1})C_4+[2k][2(k+1)])
|k,l\rangle ,
$$       $$
X_3^{(k-1,l)} X_4^{(k,l)} |k,l\rangle =
(C'_4 -(q^{2k-1}+q^{-2k+1})C_4+[2k][2(k-1)])
|k,l\rangle ,
$$
{\it where $C_4$ and $C'_4$ are the Casimir elements of
$U'_q({\rm so}_4)$ from section 2.}
\medskip

Lemma 6 can be used for evaluation of eigenvalues of Casimir elements
$C_4$ and $C'_4$ on irreducible representations $R'$ when the operators
$R(I_{21})$ and $R(I_{43})$ have eigenvalues of the nonclassical type.
Eigenvectors $| k,l\rangle$ of the operators $R'(I_{21})$ and $R'(I_{43})$
are called {\it weight vectors} of the representation $R'$. A weight
vector $| j,j'\rangle$ is called {\it a highest weight vector} if
$$
X^{(j,j')}_1 |j,j'\rangle =0,\ \ \ \ \  X^{(j,j')}_3 | j,j'\rangle =0.
$$
If $R'$ is an irreducible representation with nonclassical type eigenvalues
of the operators $R'(I_{21})$ and $R'(I_{43})$, then applying both sides
of the first and the third relations of Lemma 6 to
the vector of highest weight $|j,j'\rangle$
we derive that
$$
C_4=[ j+j'+1]_+ [ j-j']_+ I,
$$    $$
C'_4=\{ (q^{2j+1}+q^{-2j-1}) [j-j']_+ [j+j'+1]_+ -
[2j][2j+2]\} I ,
$$
where $I$ is the identity operator on the representation space.
In particular, such eigenvalues have Casimir operators of the
nonclassical type representation
$R^{\varepsilon _1,\varepsilon _2,\varepsilon _3}_{jj'}$.
\medskip

{\sl Proof of the first part of Theorem 3.} Let
$R^{\varepsilon _1,\varepsilon _2,\varepsilon _3}_{jj'}$ be a
representation of the nonclassical type on the vector space ${\cal H}$.
Then the commuting operators $R^{\varepsilon _1,\varepsilon _2,
\varepsilon _3}_{jj'} (I_{21})$ and $R^{\varepsilon _1,\varepsilon _2,
\varepsilon _3}_{jj'}(I_{43})$ are simultaneously
diagonalized. Since $q$ is not a root of unity, the
eigenvalues $(\varepsilon _1 [k+l]_+, \varepsilon _2 [k-l]_+)$
for the corresponding vector $|k,l\rangle$ are of
multiplicity 1. Let ${\cal H}'$ be a nontrivial invariant subspace
of the representation space ${\cal H}$, and let
$\sum _{k,l} \alpha _{k,l}  |k,l\rangle$ be a nonzero vector from
${\cal H}'$. Since eigenvalues
$(\varepsilon _1 [k+l]_+, \varepsilon _2 [k-l]_+)$ are of
multiplicity 1, then each $|k,l\rangle$ from this
linear combination belongs to ${\cal H}'$. Let $|k',l'\rangle$
be one of these vectors. Applying the operators $X_1$ and $X_3$
(with apropriate upper indices) to $| k',l'\rangle$ we obtain
the vector of highest weight $|j,j'\rangle$ of the
representation $R^{\varepsilon _1,\varepsilon _2,\varepsilon _3}_{jj'}$.
Applying to $|j,j'\rangle$ the operators $X_2$ and $X_4$
(with appropriate indices) we obtain all basis vectors of the
space ${\cal H}$. Hence, the representation
$R^{\varepsilon _1,\varepsilon _2,\varepsilon _3}_{jj'}$ is irreducible.
Theorem is proved.
\medskip

Now we can prove the theorem on classification of irreducible
finite dimensional representations of $U'_q({\rm so}_4)$.
\medskip

{\bf Theorem 4.} {\it If} $q$ {\it is not a root of unity, then
each irreducible finite dimensional representation $R$ of
$U'_q({\rm so}_4)$ is equivalent to one of the irreducible
representations of the classical type or to one of the
irreducible representations of the nonclassical type.}
\medskip

{\sl Proof.} Let us first prove the following assertion: {\it if
eigenvalues of the operators $R(I_{21})$ and $R(I_{43})$ of
an irreducible finite dimensional representation $R$ of
$U'_q({\rm so}_4)$ are of the classical type (that is, of the
form ${\rm i}[m]$, $m\in {\frac 12}{\Bbb Z}$),
then $R$ is equivalent to one of the
irreducible representations of the classical type}. We
diagonalize both operators $R(I_{21})$ and $R(I_{43})$ and
represent their eigenvectors in the form $|k,l\rangle$, where
$$
R(I_{21}) |k,l\rangle ={\rm i}[k+l] |k,l\rangle,\ \ \
R(I_{43}) |k,l\rangle ={\rm i} [k-l] |k,l\rangle .
$$
(These eigenvectors are called weight vectors.)
Due to Lemmas 1 and 2, there exists an eigevector of
highest weight (we denote it by $|j,j'\rangle$), that is
such that
$$
R(I_{21}) |j,j'\rangle ={\rm i}[j+j'] |j,j'\rangle,\ \ \
R(I_{43}) |j,j'\rangle ={\rm i} [j-j'] |j,j'\rangle .
$$       $$
X_1^{(j,j')}|j,j'\rangle =0,\ \ \ \ \
X_3^{(j,j')}|j,j'\rangle =0 .
$$
Applying the first and the third relations of Lemma 3 to the
vector $|j,j'\rangle$ we find eigenvalues of the Casimir
operators $R(C_4)$ and $R(C'_4)$ on the representation $R$:
$$
R(C_4) = [j+j'+1][j'-j],   \eqno (63)
$$    $$
R(C'_4) = (q^{2j+1}+q^{-2j-1})[j-j'][j+j'+1]-[2j][2j+2].   \eqno (64)
$$
Acting on the vector $|j,j'\rangle$ by the operators $X_2$
(with appropriate upper indices) we construct recursively
the vectors
$$
|j,j'-s\rangle:=X_2^{(j,j'-s+1)}\cdots X_2^{(j,j'-1)}X_2^{(j,j')}
|j,j'\rangle ,\ \ \ \ \ s=0,1,2,\cdots .
$$
Since the representation $R$ is finite dimensional and (by
Lemma 1) these vectors have different eigenvalues, there exists
smallest positive integer $n$ such that
$$
|j,j'-n\rangle:=X_2^{(j,j'-n+1)} |j,j'-n+1\rangle =0 . \eqno (65)
$$
Similarly, between vectors
$$
|j-r,j'\rangle:=X_4^{(j-r+1,j')}\cdots X_4^{(j-1,j')}X_4^{(j,j')}
|j,j'\rangle ,\ \ \ \ \ r=0,1,2,\cdots ,
$$
there exists a nonzero vector with smallest positive integer
$m$ such that
$$
|j-m,j'\rangle:=X_4^{(j-m+1,j')} |j-m+1,j'\rangle =0 . \eqno (66)
$$
Then using the second and the fourth relations of Lemma 3, we
find from (63)--(66) that
$$
X_1^{(j,j'-n)} X_2^{(j,j'-n+1)} |j,j'-n+1\rangle
$$           $$
=[n][n-2j'-1](q^{j'-j-n}+q^{-j'+j+n})(q^{j+j'-n+1}+q^{-j-j'+n-1})
|j,j'-n+1\rangle = 0 ,
$$     $$
X_3^{(j-m,j')} X_4^{(j-m+1,j')} |j-m+1,j'\rangle
$$           $$
=[m][m-2j-1](q^{m-j-j'-1}+q^{-m+j+j'+1})(q^{m-j+j'}+q^{-m+j-j'})
|j-m+1,j'\rangle = 0 .
$$
Therefore, $[n-2j'-1]=0$ and $[m-2j-1]=0$, that is
$$
m=2j+1, \ \ \ \ \ n=2j'+1 .
$$

Now we act successively on the vectors $|j,j'\rangle , |j,j'-1\rangle ,
|j,j'-2\rangle ,\cdots ,|j,-j'\rangle $ by the operators $X_4$ with
the appropriate upper indices. As a result, we construct the
vectors
$$
|j,j'\rangle ,\ \ \  \ |j-1,j'\rangle ,\ \ \  \cdots \ \ \
|-j,j'\rangle ,
$$          $$
|j,j'-1\rangle , \ \  \ |j-1,j'-1\rangle , \ \  \cdots \ \
|-j,j'-1\rangle ,
$$          $$
\cdots \ \ \ \ \ \ \ \  \cdots \ \ \ \ \ \ \ \ \cdots \ \ \ \ \ \ \ \
\cdots
$$        $$
|j,-j'\rangle ,\ \ \  \ |j-1,-j'\rangle ,\ \ \  \cdots \ \ \
|-j,-j'\rangle
$$
for which
$$
R(I_{21})|k,l\rangle = {\rm i}[k+l] |k,l\rangle ,\ \ \ \  \ \
R(I_{43})|k,l\rangle = {\rm i}[k-l] |k,l\rangle .
$$
We can find how the operators $X_i$, $i=1,2,3,4$,
with appropriate indices act on these vectors:
$$
X_2^{(k,l)} |k,l\rangle =|k,l-1\rangle ,\ \ \ \ \ \
X_4^{(k,l)} |k,l\rangle = |k-1,l\rangle ,
$$       $$
X_1^{(k,l)} |k,l\rangle =X_1^{(k,l)}X_2^{(k,l+1)} |k,l+1\rangle
=(R(C'_4)-(q^{(2l+1}+q^{-2l-1})R(C_4)+[2l+2][2l])|k,l+1\rangle ,
$$         $$
X_3^{(k,l)} |k,l\rangle =X_3^{(k,l)}X_4^{(k+1,l)} |k+1,l\rangle
=(R(C'_4)+(q^{(2k+1}+q^{-2k-1})R(C_4)+[2k+2][2k])|k+1,l\rangle .
$$
Putting here the explicit expression for the Casimir operators,
substituting these expressions for $X_i^{(k,l)} |k,l\rangle$
to (55)--(58) and considering (55)--(58) as a system of
linear equations with unknown $R(I_{32}) |k,l\rangle$,
$R(I_{31}) |k,l\rangle$, $R(I_{42}) |k,l\rangle$, $R(I_{41}) |k,l\rangle$
we solve this system and find that
$$
R(I_{32})|k,l\rangle =\frac{1}{(q^{k+l}+q^{-k-l})(q^{k-l}+q^{-k+l})}
\{ |k-1,l\rangle +|k,l-1\rangle -
$$      $$
-(q^{j-l}+q^{j+l})(q^{j+l+1}+q^{-j-l-1}) [j'-l][j'+l+1]
|k,l+1\rangle -
$$        $$
-(q^{j'-k}+q^{j+k})(q^{j'+k+1}+q^{-j'-k-1}) [j-k][j+k+1]
|k+1,l\rangle \} .
$$
Thus, the vectors $|k,l\rangle$, $-j\le k\le j$, $-j'\le l\le j'$,
constitute a basis of the representation space ${\cal H}$.
Introducing a new basis $\{ |k,l\rangle '\}$ such that
$$
|k,l\rangle =(-1)^{k+l} \prod _{r=-j}^k (q^{j'+r}+q^{-j'-r})^{-1}[j+r]^{-1}
\prod _{s=-j'}^k (q^{j+s}+q^{-j-s})^{-1}[j'+s]^{-1}
|k,l\rangle '
$$
we shall obtain for $R(I_{21})$, $R(I_{32})$, $R(I_{43})$
the operators of the irreducible representation
$R_{jj'}$ of the classical type from section 4. Thus, in
this case Theorem is proved.

Now we prove the second part of the theorem which can be formulated
as follows: {\it if eigenvalues of the operators $R(I_{21})$ and
$R(I_{43})$ of an irreducible finite dimensional representation $R$
are of the nonclassical type, that is of the form $\pm [m]_+$,
$m\in {\frac 12}{\Bbb Z}$, $m\not\in {\Bbb Z}$, then
$R$ is equivalent to one of the irreducible representations of the
nonclassical type.}

We first prove that if eigenvalues of the operators $R(I_{21})$ and
$R(I_{43})$ are of the form $[m]_+$ (only sign + is taken), then
$R$ is equivalent to one of the irreducible representations of the
nonclassical type. A proof is similar to that of the previous case
and for this reason we do not give details.

Due to Lemmas 4 and 5, there exists an eigenvector of highest
weight $|j,j'\rangle$ such that
$$
R(I_{21}) |j,j'\rangle =[j+j']_+ |j,j'\rangle,\ \ \
R(I_{43}) |j,j'\rangle =[j-j']_+ |j,j'\rangle .
$$       $$
X_1^{(j,j')}|j,j'\rangle =
X_3^{(j,j')}|j,j'\rangle =0 .
$$
Applying relations of Lemma 6 to $|j,j'\rangle$ we find
eigenvalues of the Casimir operators $R(C_4)$ and $R(C'_4)$:
$$
R(C_4) = [j+j'+1 ]_+ [ j-j' ]_+ ,
$$    $$
R(C'_4) = (q^{2j+1}+q^{-2j-1})[ j-j']_+ [ j+j'+1]_+ -[2j][2j+2].
$$
Now we construct recursively the vectors
$$
|j-r,j'-s\rangle:=X_4^{(j-r+1,j'-s)} \cdots X_4^{(j,j'-s)}
X_2^{(j,j'-s+1)}\cdots X_2^{(j,j'-1)}X_2^{(j,j')}
|j,j'\rangle , \eqno (67)
$$      $$
r,s=0,1,2,\cdots .
$$
By Lemma 4, for these vectors we have
$$
R(I_{21}) |k,l\rangle =[k+l]_+ |k,l\rangle ,\ \ \ \
R(I_{43}) |k,l\rangle =[k-l]_+ |k,l\rangle  . \eqno (68)
$$
These vectors satisfy the relations
$$
X_2^{(k,l)} |k,l\rangle =|k,l-1\rangle ,\ \ \ \ \
X_4^{(k,l)} |k,l\rangle = |k-1,l\rangle , \eqno (69)
$$       $$
X_1^{(k,l)} |k,l\rangle =X_1^{(k,l)}X_2^{(k,l+1)} |k,l+1\rangle
$$        $$
=(R(C'_4)-(q^{(2l+1}+q^{-2l-1})R(C_4)+[2l+2][2l])|k,l+1\rangle ,
\eqno (70)
$$         $$
X_3^{(k,l)} |k,l\rangle =X_3^{(k,l)}X_4^{(k+1,l)} |k+1,l\rangle
$$        $$
=(R(C'_4)-(q^{(2k+1}+q^{-2k-1})R(C_4)+[2k+2][2k])|k+1,l\rangle .
\eqno (71)
$$

Since the operators $X_1^{(k,l)}$, $X_2^{(k,l)}$, $X_3^{(k,l)}$,
$X_4^{(k,l)}$ acting on the vector $|k,l\rangle$ determine the
action of the operators  $R(I_{32})$, $R(I_{31})$, $R(I_{42})$,
$R(I_{41})$ on this vector, then the vectors
$$
|j,j'\rangle ,\ \ \ \ \ |j-1,j'\rangle ,\ \ \ \ \ |j-2,j'\rangle ,
\ \ \ \ \ \cdots ,
$$      $$
|j,j'-1\rangle , \ \ \ \ |j-1,j'-1\rangle , \ \ \ \ |j-2,j'-1\rangle ,
 \ \ \ \ \cdots , \eqno (72)
$$      $$
|j,j'-2\rangle , \ \ \ \ |j-1,j'-2\rangle , \ \ \ \ |j-2,j'-2\rangle ,
 \ \ \ \ \cdots ,
$$    $$
\cdots \ \ \ \ \ \ \ \cdots \ \ \ \ \ \ \ \cdots \ \
\ \ \ \ \ \cdots
$$
span an invariant subspace in the representation space ${\cal H}$.
Since the representation $R$ is irreducible, they span the whole
space ${\cal H}$. It follows from (68) that only pairs of vectors
$|k,l\rangle$ and $|-k,-l\rangle$ have the same eigenvalues for
the operators $R(I_{21})$ and $R(I_{43})$.

In order to determine which possibilities exist for the representation
$R$, we make as follows. For definiteness we suppose that $j$ is
half-integral. Using formula (67) we first create the set of all
possible vectors $|k,l\rangle$, which does not contains pairs
$|k,l\rangle$ and $|-k,-l\rangle$. For example, we create all the
vectors $|k,l\rangle$ with $k>0$. This set contains the vector
$| \frac 12 ,0\rangle$. There are two possibities:
\medskip

(a) The vectors $| \frac 12 ,0\rangle$  and $| -\frac 12 ,0\rangle$
are linearly dependent, that is,
$X_4^{(1/2,0)}| \frac 12 ,0\rangle =a|\frac 12 ,0\rangle$;

(b) the vectors $| \frac 12 ,0\rangle$  and $| -\frac 12 ,0\rangle$
are linearly independent.
\medskip

\noindent
In the case (a) all pairs $|k,l\rangle$ and $|-k,-l\rangle$
consist of linearly dependent vectors. The reason of this is
that
$$
X^{(k,l)}_1= X^{(-k,-l)}_2,\ \ \ \ \
X^{(k,l)}_3=X^{(-k,-l)}_4          \eqno (73)
$$
(as it follows from expressions (59)--(62) for the operators
$X_i$, $i=1,2,3,4$). Therefore, for every positive integrals $r$ and
$s$ the vectors
$$
X^{(r+1/2,s)}_3 \cdots X^{(r+1/2,0)}_3
X^{(r-1/2,0)}_1 \cdots X^{(1/2,0)}_1 |1/2,0\rangle ,
$$          $$
X^{(-r-1/2,-s)}_4 \cdots X^{(-r-1/2,0)}_4
X^{(-r+1/2,0)}_2 \cdots X^{(-1/2,0)}_2 |-1/2,0\rangle
$$
are linear dependent with the same constant $a$.
The constant $a$ can be explicitly calculated. Namely,
since $X_4^{(1/2,0)} |1/2,0\rangle =| -1/2,0\rangle =a| 1/2,0\rangle$
and $X_3^{(-1/2,0)}=X_4^{(1/2,0)}$, we have
$$
X_3^{(-1/2,0)} |-1/2,0\rangle =X_3^{(-1/2,0)}X_4^{(1/2,0)} |1/2,0\rangle
=R(C'_4 -2C_4-1) |1/2,0\rangle =
$$    $$
=X_3^{(-1/2,0)} a|1/2,0\rangle =a X_4^{(1/2,0)} |1/2,0\rangle =
a^2 |1/2,0\rangle ,
$$
that is, $a^2=R(C'_4-2C_4-1)$. This means that $a$ is determined up to a
sign. Using expressions for values of the Casimir elements of
$U'_q({\rm so}_4)$ on irreducible representations with highest weight
vector $| j,j'\rangle$, we find that
$$
a=\varepsilon _3 (q-q^{-1}) [j+1/2][j'+1/2],
$$
where $\varepsilon _3$ takes one of the values $\pm 1$.

Thus, we recieved the following set of linear independent
vectors of the representation space ${\cal H}$:
$$
|j,j'\rangle ,\ \ \ \ \ |j-1,j'\rangle ,\ \ \ \ \ \cdots ,
\ \ \ \ \  | 1/2 ,j'\rangle ,
$$      $$
|j,j'-1\rangle , \ \ \ \ |j-1,j'-1\rangle , \ \ \ \ \cdots ,
 \ \ \ \ | 1/2 ,j'-1\rangle  ,
$$      $$
\cdots \ \ \ \ \ \ \ \ \ \cdots \ \ \ \ \ \ \ \ \ \cdots \ \ \ \
\ \ \ \ \ \cdots
$$      $$
|j,-j'\rangle , \ \ \ \ |j-1,-j'\rangle , \ \ \ \ \cdots ,
 \ \ \ \ | 1/2 ,-j'\rangle  .
$$
These vectors constitute a basis  of the space ${\cal H}$. We
introduce a new basis $\{ |k,l\rangle '\}$ such that
$$
|k,l\rangle = \prod _{r=1/2}^k (q^{j'+r}-q^{-j'-r})^{-1}[j+r]^{-1}
\prod _{s=-j'}^k (q^{j+s}-q^{-j-s})^{-1}[j'+s]^{-1}
|k,l\rangle ' .
$$
Rewriting the relations (69)--(71) for this new basis we obtain
a system of linear equations with unknown
$R(I_{32})|k,l\rangle '$, $R(I_{42})|k,l\rangle '$,
$R(I_{41})|k,l\rangle '$, $R(I_{31})|k,l\rangle '$. Solving
this system we find that the operator $R(I_{32})$ acts on the
basis vectors by the formulas for the irreducible
nonclassical type representation $R_{jj'}^{(+,+,+)}$ or the
irreducible nonclassical type representation
$R_{jj'}^{(+,+,-)}$.

Now we consider case (b). Due to the relation (73) we conclude
that the vectors
$$
|k,l\rangle ,\ \ \ \ \ k=j,j-1,j-2,\cdots ,-j,\ \ \ \
l=j',j'-1'j'-2,\cdots ,-j',
$$
are linearly independent and constitute a basis of the
representation space ${\cal H}$. Solving the system of
linear equations (69)--(71) we obtain a representation of
$U'_q({\rm so}_4)$ equivalent to one of the representations
$R^{(\pm {\rm i},\pm 1)}_{jj'}$ or
$R^{(\pm 1,\pm {\rm i})}_{jj'}$ from section 6. These
representations are reducible. So, case (b) is not possible
for our irreducible representation $R$.

We have considered the case when all eigenvalues of the
operators $R(I_{21})$ and $R(I_{43})$ are of the form
$[m]_+$ (with sign +). However, due to automorphisms
$\psi _1$ and $\psi _2$, mapping $I_{21}\to -I_{21}$ and
$I_{43}\to -I_{43}$, respectively, and concerving all other
generating elements in $\{ I_{21},I_{32},I_{43} \}$, to every such
irreducible representation $R$ there correspond the
representations $R^{(+,-)}=R\circ \psi _2$,
$R^{(-,+)}=R\circ \psi _1$, $R^{(-,-)}=R\circ \psi _1\psi _2$
such that
$$
R^{(+,-)}(I_{21})|k,l\rangle =[k+l]_+ |k,l\rangle ,\ \ \ \
R^{(+,-)}(I_{43})|k,l\rangle =-[k-l]_+ |k,l\rangle , \eqno (74)
$$         $$
R^{(-,+)}(I_{21})|k,l\rangle =-[k+l]_+ |k,l\rangle ,\ \ \ \
R^{(-,+)}(I_{43})|k,l\rangle =[k-l]_+ |k,l\rangle , \eqno (75)
$$         $$
R^{(-,-)}(I_{21})|k,l\rangle =-[k+l]_+ |k,l\rangle ,\ \ \ \
R^{(-,-)}(I_{43})|k,l\rangle =-[k-l]_+ |k,l\rangle . \eqno (76)
$$
Conversely, to any of the representations
$R^{(+,-)}$, $R^{(-,+)}$, $R^{(-,-)}$ with these properties
there corresponds a unique representation $R$ such that
$$
R(I_{21})|k,l\rangle =[k+l]_+ |k,l\rangle ,\ \ \ \
R(I_{43})|k,l\rangle =[k-l]_+ |k,l\rangle  \eqno (77)
$$
for all eigenvactors $|k,l\rangle$. This means that the
classification of irreducible representations $R$ with
property (77) automatically leads to the classification of
irreducible representations with any of the properties
(74)--(76) and vise versa. Therefore, any of irreducible
representations of $U'_q({\rm so}_4)$ with one of the
properties (74)--(76) is equivalent to one of the
irreducible representations of the nonclassical type.
Theorem is proved.
\medskip

In section 4 we constructed the homomorphism $\phi :
U'_q({\rm so}_4)\to
{\hat U}_q({\rm sl}_2)^{\otimes 2,{\rm ext}}$ (see Theorem 2).
Now we are able to prove more strong assertion:
\medskip

{\bf Corollary.} {\it If} $q$ {\it is not a root of unity,
then the homomorphism $U'_q({\rm so}_4)\to
{\hat U}_q({\rm sl}_2)^{\otimes 2,{\rm ext}}$ of
Theorem 2 is injective.}
\medskip

{\sl Proof.} If the assertion of Corollary is not true,
then there exists nonzero element $a\in U'_q({\rm so}_4)$
such that $\phi (a)=0$. Then for any finite dimensional
representation $T$ of the algebra
${\hat U}_q({\rm sl}_2)^{\otimes 2,{\rm ext}}$  we
have $T(\phi (a))=0$. Taking the representations
$T_{jj'}^{(\pm ,\pm )}$, $T_{jj'}^{(\pm {\rm i},\pm )}$
and $T_{jj'}^{(\pm ,\pm {\rm i})}$,
where $T^{(\pm \varepsilon _1 ,\pm \varepsilon _2)}_{jj'}\equiv
T^{(\pm \varepsilon _1)}_j \otimes T^{(\pm \varepsilon _2)}_{j'}$
(see sections 6 and 7), as a representation $T$ we obtain
$$
T_{jj'}^{(\pm ,\pm )}(\phi (a))=0,\ \ \ \
T_{jj'}^{(\pm {\rm i},\pm )} (\phi (a))=0,\ \ \ \
T_{jj'}^{(\pm ,\pm {\rm i})} (\phi (a))=0
$$
for all admissable values of $j$ and $j'$. As we have
seen above, any irreducible finite dimensional
representation of $U'_q({\rm so}_4)$ is equivalent
to the representation $T_{jj'}^{(1,1)}\circ \phi$
with appropriate values of $j$ and $j'$ or to one of
irreducible constituents of the representations
$T_{jj'}^{(\pm {\rm i},\pm )} \circ \phi $ and
$T_{jj'}^{(\pm ,\pm {\rm i})} \circ \phi $. This
means that $R(a)=0$ for any irreducible finite
dimensional representation of $U'_q({\rm so}_4)$.
But it was shown in [13] (see also [2]) that irreducible finite
dimensional representations of $U'_q({\rm so}_4)$
separate elements of this algebra. Thus, for our
element $a$ these exists an irreducible finite
dimensional representation $R$ such that $R(a)\ne 0$.
This contradiction proves Corollary.
\bigskip

\noindent
{\sf 9. COMPLETE REDUCIBILITY OF FINITE DIMENSIONAL
REPRESENTATIONS}
\medskip

It was proved in [29] that if $q$ is not a root of unity,
then every finite dimensional
representation of the algebra $U'_q({\rm so}_3)$ is
completely reducible. The aim of this section is to
prove the corresponding theorem for the algebra
$U'_q({\rm so}_4)$.
\medskip

{\bf Theorem 5.} {\it If} $q$ {\it is not a root of unity, then each
finite dimensional representation of $U'_q({\rm so}_4)$ is
completely reducible.}
\medskip

\noindent
{\sl Proof.}
In order to prove this theorem it is enough to show that every finite
dimensional representation $R$ of $U'_q({\rm so}_4)$, containing only
two irreducible constituents, is completely reducible.

There are three possibilities for two irreducible constituents:
\medskip

(a) both representations are of the classical type;

(b) both representations are of the nonclassical type;

(c) representations belong to different types.
\medskip

\noindent
Each case will be proved separately.
\medskip

{\bf Case (a).}
We shall use in the proof the following properties of the operators
$X_i^{(k,l)}$ from (55)--(58), which are derived from Lemmas 1--3:
\medskip

(A) Let $|k,l\rangle $ be such as in Lemma 1. Acting on $|k,l\rangle$
by the operators $X_i$, $i=1,2,3,4$, with appropriate upper indices,
we can obtain a vector $|k,l\rangle '$ with eigenvalues of the operators
$R(I_{21})$ and $R(I_{43})$ coinciding with those of the vector
$|k,l\rangle $. Then $|k,l\rangle '=a|k,l\rangle $ for some complex number
$a$.

(B) Let $|k,l\rangle$ be such as in Lemma 1. If $X_1^{(k,l)}|k,l\rangle
=|k,l+1\rangle =0$ and $|k',l\rangle$ is another weight vector of the
operators $R(I_{21})$ and $R(I_{43})$ obtained by action of the operators
$X_i$, $i=1,2,3,4$, with appropriate upper indices, then
$X_1^{(k',l)}|k',l\rangle =|k',l+1\rangle =0$.
The same assertion is valid for the vectors
$X_2^{(k,l)}|k,l\rangle$, $X_3^{(k,l)}|k,l\rangle$ and
$X_4^{(k,l)}|k,l\rangle$. This means that by acting on $|k,l\rangle $
by the operators $X_i$, $i=1,2,3,4$, with appropriate upper indices, we
obtain the set of nonzero vectors $|k',l'\rangle $ such that their
values $(k',l')$ constitute a parallelogram.
\medskip

First let us consider the subcase when two constituents are
equivalent. We denote them by $R_{jj'}$ and $R'_{jj'}$.
Since restriction of the representation $R$ to the subalgebra
$U'_q({\rm so}_3)$ is completely reducible,
there exists a basis in the space of the
representation $R$ consisting of eigenvectors for both operators
$R(I_{21})$ and $R(I_{43})$. In this basis there exist exactly two
vectors of highest weight. Let
$| j,j'\rangle $ and  $|j,j'\rangle '$ be these vectors.
We create two sets of vectors
$$
X^r_2 X^t_4 |j,j'\rangle ,\ \ \ r,t=0,1,2,\cdots ,\ \ \ \
{\rm and}\ \ \ \
X^r_2X^t_4 |j,j'\rangle ', \ \ \ r,t=0,1,2,\cdots ,
$$
where the operators $X_2$ and $X_4$ are taken with the appropriate upper
indices. Due to the properties of the operators $X_2$ and $X_4$, these
two sets span two subspaces $V_1$ and $V_2$ which are
invariant with respect to
the operators $R(I_{21})$, $R_(I_{32})$ and $R(I_{43})$. Moreover, we
have $V=V_1\oplus V_2$ and the theorem is proved in this case.

Now suppose that two constituents (denote them by
$R_{jj'}$ and $R_{ss'}$) of
our reducible representation $R$ are not equivalent and that $R_{jj'}$
is realized in an invariant subspace $V_1$. As in the previous subcase,
the operators $R(I_{21})$ and $R(I_{43})$ can be simultaneously
diagonalized.
We represent the whole
representation space $V$ in the form $V=V_1\oplus V_2$, where $V_2$ have
a basis consisting of eigenvectors of the operators $R(I_{21})$ and
$R(I_{43})$. Two subcases are possible in this case:
\medskip

(I) highest weight of the representation $R_{ss'}$ is not a weight of
$R_{jj'}$;

(II) highest weight of $R_{ss'}$ is a weight of $R_{jj'}$.
\medskip

\noindent
In the first subcase let $|s ,s' \rangle '$ be a vector of
highest weight for the representation $R_{ss'}$ in the subspace
$V_2$ (its eigenvalue is of multiplicity 1). We create the set of
vectors
$$
X^r_2X^t_4 |s ,s' \rangle ',\ \ \ r,t=0,1,2,\cdots ,
$$
where $X_2$ and $X_4$ are taken with the appropriate upper indices.
Then these vectors span a subspace invariant with respect to action of
the operators $X_i$, $i=1,2,3,4$, taken with appropriate indices.
Therefore, this subspace is invariant with respect to the
operators $R(I_{21})$, $R(I_{32})$, $R(I_{43})$, and no of these
basis vectors belong to $V_1$. This means that the representation $R$ is
completely reducible in this subcase.

In the subcase (II) we consider the eigenspace $V_{k,l}$ of the
operators $R(I_{21})$ and $R(I_{43})$ with eigenvalues $k+l$ and
$k-l$, respectively,
such that
$$
{\rm dim}\, V_{k,l}=2,\ \ \ \
{\rm dim}\, V_{k,l+1}< 2,\ \ \ \
{\rm dim}\, V_{k+1,l}< 2.
$$
(This means that the highest weight of the representation $R_{ss'}$ is
$(k+l,k-l)$,that is, $k=s$, $l=s'$.) Then ${\rm dim}\, V_{k,l+1}=1$
or/and ${\rm dim}\,
V_{k+1,l}=1$. Let ${\rm dim}\, V_{k+1,l}=1$ (if ${\rm dim}\,
V_{k,l+1}=1$, then a proof is the same). Then in $V_{k,l}$ there exists a
vector $|k,l\rangle '$ such that $X_3^{(k,l)} |k,l\rangle '=0$.
Let us show that $X_1^{(k,l)} |k,l\rangle '=0$. Suppose that
$X_1^{(k,l)}|k,l\rangle '=\alpha |k,l+1\rangle '$,
$\alpha \ne 0$. Then
we act successively on the vector $|k,l+1\rangle '$ by the operators
$X_i$, $i=1,2,3,4$, with the appropriate upper indices and obtain the set
of vectors $|k',l'\rangle '$ which span invariant subspace $V'$. Due
to the property (B) of the operators $X_i$, this set of vectors does not
contain vectors $|k+1,l'\rangle '$ with some values of $l'$.
This means that $V'$ is a nontrivial
subspace of the representation space. A highest weight vector of the
subspace $V'$ does not coincide with $(s,s')$ and $(j,j')$. Thus,
the whole representation $R$ contains the third irreducible constituent.
It is a contradiction which show that
$X_1^{(k,l)} |k,l\rangle '=0$.

We constructed the vector $| k,l\rangle '$ such that
$X_1^{(k,l)} |k,l\rangle '=0$ and $X_3^{(k,l)} |k,l\rangle '=0$.
Let $V_2$ be a subspace of the representation
space $V$ spanned by the vectors
$$
X_2^rX_4^t |k,l\rangle ' ,\ \ \ \ r,t=0,1,2,\cdots .
$$
According to the properties of the operators $X_2$ and $X_4$,
then $V_2$ is invariant subspace for the operators $R(I_{21})$,
$R(I_{32})$, $R(I_{43})$ and $V=V_1\oplus V_2$, where $V_1$ is the
subspace of the irreducible representation $R_{jj'}$. Theorem is proved
in this subcase.
\medskip

\noindent {\bf Case (b).}
We distinguish here the following subcases:
\medskip

(I) Two constituents are equivalent;

(II) two constituents are of the form $R_{jj'}^{(\varepsilon _1,
\varepsilon _2,\varepsilon _3)}$ and $R_{jj'}^{(\varepsilon _1',
\varepsilon _2', \varepsilon _3')}$, where $\varepsilon _1\ne
\varepsilon _1'$ or/and $\varepsilon _2\ne \varepsilon _2'$;

(III) two constituents are of the form $R_{jj'}^{(\varepsilon _1,
\varepsilon _2,\varepsilon _3)}$ and $R_{jj'}^{(\varepsilon _1,\varepsilon
_2,
-\varepsilon _3)}$;

(IV) two constituents are $R_{jj'}^{(\varepsilon _1,
\varepsilon _2,\varepsilon _3)}$ and $R_{ss'}^{(\varepsilon _1',\varepsilon
_2',
\varepsilon _3')}$, $(j,j')\ne (s,s')$.
\medskip

\noindent
For subcase (I) a proof is such as in the first part of case (a).
Let us consider subcase (II). Let $\varepsilon _1\ne \varepsilon _1'$.
In the representation space $V$, there exist two
linearly independent vectors $| j,j'\rangle $
and $|j,j'\rangle '$ which are of highest weights, that is,
$$
X_1^{(j,j')}|j,j'\rangle  = X_1^{(j,j')}|j,j'\rangle ' = 0,\ \ \ \
X_3^{(j,j')}|j,j'\rangle  = X_3^{(j,j')}|j,j'\rangle ' = 0,
$$
and such that
$$
R(I_{21})|j,j'\rangle  =\varepsilon _1 [j+j']_+|j,j'\rangle ,\ \ \ \
R(I_{21})|j,j'\rangle ' =\varepsilon '_1 [j+j']_+|j,j'\rangle .
$$
We create two sets of vectors
$$
X_2^rX_4^t |j,j'\rangle ,\ \ \ \ \ r,t=0,1,2,\cdots ,\ \ \ \
{\rm and}\ \ \ \
X_2^rX_4^t |j,j'\rangle ',\ \ \ \ \ r,t=0,1,2,\cdots .
$$
No nonzero vector of the first set is multiple of some vector of the
second set (since otherwise these two sets span the same vector
subspace of the representation space $V$). These two sets of vectors
span two invariant linear subspaces $V_1$ and $V_2$
of $V$. Since $V=V_1\oplus V_2$, then the
representation $R$ is a direct sum of the subrepresentations
$R_{jj'}^{(\varepsilon _1,\varepsilon _2,\varepsilon _3)}$ and
$R_{jj'}^{(\varepsilon _1',\varepsilon _2',
\varepsilon _3')}$.

Let us prove the theorem in subcase (III). We suppose that the
representation $R_{jj'}^{(\varepsilon _1,\varepsilon _2,\varepsilon _3)}$
is realized on an invariant subspace $V_1$ and denote weight vectors for
this representation by $|k,l\rangle $.  We have
$R(I_{32})|1/2,0\rangle =  \varepsilon _3a |1/2,0\rangle +\cdots $,
where $a$ is the appropriate constant and
a linear combination of weight vectors with
weights different from that of the vector $|1/2,0\rangle $ is denoted by
dots. In the space $V$ of the representation $R$ there exists another
vector $|1/2,0\rangle '$ such that
$$
R(I_{32})|1/2,0\rangle '= - \varepsilon _3a |1/2,0\rangle ' +
r|1/2,0\rangle  +\cdots ,
$$
where dots mean the same as above. Then we easily verify that
$$
R(I_{32})\left( |1/2,0\rangle '
-\frac r{2a\varepsilon _3} |1/2,0\rangle  \right) = -a \varepsilon _3
\left( |1/2,0\rangle ' -\frac r{2a\varepsilon _3} |1/2,0\rangle  \right)
+\cdots
$$
with the same meaning for dots. Denoting the vector
$|1/2,0\rangle ' -\frac r{2a\varepsilon  _3} |1/2,0\rangle $ by
$|1/2,0\rangle ''$ we create the vectors
$$
X_1^rX_3^t |1/2,0\rangle '',\ \ \ \ \ r,t=0,1,2,\cdots ,
$$
then take the vector of highest weight in this set (we denote it
by $|j,j'\rangle ''$) and create the vectors
$$
X_1^rX_3^t |j,j'\rangle '',\ \ \ \ \ r,t=0,1,2,\cdots .
$$
The linear subspace $V_2$ spanned by the last vectors is invariant and
irreducible. Since $V=V_1\oplus V_2$, the theorem is proved in
this subcase.

The subcase (IV) is proved in the same way as the second part of the
previous case.
\medskip

\noindent
{\bf Case (c).} First we note that we can diagonalize simultaneously
the operators $T(I_{21})$ and $T(I_{43})$.
Let $V=V_1\oplus V_2$ be the decomposition of the representation space
into the direct sum of subspaces $V_1$ and $V_2$ spanned by eigenvectors
of the classical types (with eigenvalues of the type $[k+l]$ and
$[k-l]$) and of the nonclassical type (with eigenvalues of the type
$[l+k]_+$ and $[k-l]_+$), respectively. Now we take a vector of highest
weight in the subspace $V_1$ (denote it by $|j,j'\rangle $) and a
vector of highest weight
in the subspace $V_2$ (denote it by $|s,s'\rangle '$) and then create two
sets of vectors
$$
X_2^rX_4^t |j,j'\rangle ,\ \ \ \ \ r,t=0,1,2,\cdots ,\ \ \ \
{\rm and}\ \ \ \
X_2^rX_4^t |s,s'\rangle ',\ \ \ \ \ r,t=0,1,2,\cdots .
$$
In fact, they span the vector subspaces $V_1$ and $V_2$, respectively.
Using properties of the operators $X_2$ and $X_4$, we conclude that
the subspaces $V_1$ and $V_2$ are invariant with respect to the operators
$R(I_{21})$, $R(I_{32})$ and $R(I_{43})$. Theorem is proved.
\bigskip

\noindent
{\sf 10. TENSOR PRODUCTS OF REPRESENTATIONS}
\medskip

As in the case of the algebra $U'_q({\rm so}_3)$ (see [21]), the
homomorphism of Theorem 2 allows us to determine tensor products of
irreducible finite dimensional representations of the algebra
$U'_q({\rm so}_4)$ and decompose them into irreducible constituents.

Let us explain this on the example of the classical type irreducible
representations of $U'_q({\rm so}_4)$. Let $T^{(1,1)}_{jj'}=T^{(1)}_j
\otimes T^{(1)}_{j'}$ be the irreducible representation of
$U_q({\rm sl}_2)^{\otimes 2,{\rm ext}}$. Then the tensor product
$T^{(1,1)}_{jj'} \otimes T^{(1,1)}_{ss'}$ is well defined representation
of the algebra $U_q({\rm sl}_2)^{\otimes 2}$. Thus, we have to
determine the operators $(T^{(1,1)}_{jj'} \otimes T^{(1,1)}_{ss'})(x_i)$,
$i=1,2$. For this we use the determined operators
$(T^{(1,1)}_{jj'} \otimes T^{(1,1)}_{ss'})(c_i)$, $i=1,2$, where $c_i$
are Casimir elements from section 3. We define the
operators $(T^{(1,1)}_{jj'} \otimes T^{(1,1)}_{ss'})(x_i)$ as solutions of
equations
$$
q^{-1}\{ (T^{(1,1)}_{jj'} \otimes T^{(1,1)}_{ss'}) (x_i)\} ^4 -
 (T^{(1,1)}_{jj'} \otimes T^{(1,1)}_{ss'})(c_i) (q-q^{-1})^2
\{ (T^{(1,1)}_{jj'} \otimes T^{(1,1)}_{ss'}) (x_i)\} ^2 +q=0
\eqno (78)
$$
(see the equation for the elements $x_i$ in section 3).
In order to find these solutions we may diagonalize the operators
$( T^{(1,1)}_{jj'} \otimes T^{(1,1)}_{ss'}) (c_i)$. Then solutions of
equations (78) can be easily calculated. Composing the representation
$T^{(1,1)}_{jj'} \otimes T^{(1,1)}_{ss'}$ of
$U_q({\rm sl}_2)^{\otimes 2,{\rm ext}}$ with the homomorphism $\phi$
from Theorem 2 we obtain the representation
$$
R_{jj'}\otimes R_{ss'} =\{
T^{(1,1)}_{jj'} \otimes T^{(1,1)}_{ss'} \} \circ \phi
$$
which is treated as the tensor product of irreducible representations
$R_{jj'}$ and $R_{ss'}$ of the algebra $U'_q({\rm so}_4)$.

For the representation $T^{(1,1)}_{jj'} \otimes T^{(1,1)}_{ss'} $ of
$U_q({\rm sl}_2)^{\otimes 2,{\rm ext}}$ we have the decomposition
$$
T^{(1,1)}_{jj'} \otimes T^{(1,1)}_{ss'} =\sum _{k=|j-s|}^{j+s}
\sum _{k'=|j'-s'|}^{j'+s'} \oplus T^{(1,1)}_{kk'} .
$$
Composing the left and the right hand sides of this relation with the
homomorphism $\phi$ from Theorem 2:
$$
(T^{(1,1)}_{jj'} \otimes T^{(1,1)}_{ss'})\circ \phi =\left( \sum _k
\sum _{k'}\oplus T^{(1,1)}_{kk'} \right) \circ \phi
$$
we obtain the decomposition of the tensor product
$R_{jj'}\otimes R_{ss'}$:
$$
R_{jj'}\otimes R_{ss'} = \sum _{k=|j-s|}^{j+s} \sum _{k'=|j'-s'|}^{j'+s'}
\oplus
R_{kk'} .
$$
As in the case of the algebra $U'_q({\rm so}_3)$ in [21],
we can also determine in the similar way tensor products of irreducible
representations of the classical and the nonclassical types and tensor
products of irreducible representations of the nonclassical type.
\bigskip

\noindent
{\sf ACKNOWLEDGMENTS}

One of the authors (AUK) would like to thank International E. Schr\"odinger
Institute for Mathematical Physics in Vienna and Department of Mathematics
of FNSPE, Czech Technical University, for the hospitality. Part of this
paper was included into the report of this author in the Program
"Algebraic Groups, Theory of Invariants, and Applications" in
International E. Schr\"odinger Institute for Mathematical Physics.

\end{document}